\documentclass[authoryear]{sn-jnl}

\usepackage[english]{babel}
\usepackage[utf8]{inputenc}

\usepackage{lipsum, mathrsfs}
\usepackage{
    mathtools,
    amsmath,
    amssymb,
    amsthm,
    enumerate,
    mdframed,
    braket,
    physics,
    optidef,
    algorithm,
    booktabs,
    hyperref,
    natbib
}
\usepackage[noend]{algpseudocode}
\usepackage{dsfont}
\usepackage{cleveref}
\usepackage{lucassymbs, lucasdelims}

\begin{document}

\title{Robust Parallel Machine Scheduling: A Comparison of Robustness Measures and Local Search Approaches}

\author*[1]{\fnm{Casper} \sur{Loman}}\email{c.loman@uu.nl}

\author[1]{\fnm{Loriana} \sur{Pascual}}

\author[1]{\fnm{Bram} \sur{Elderhorst}}

\author[1]{\fnm{Marjan} \sur{van den Akker}}\email{j.m.vandenakker@uu.nl}

\author[2]{\fnm{Roel} \sur{van den Broek}}\email{roel.vandenbroek@ns.nl}

\author[1,2]{\fnm{Han} \sur{Hoogeveen}}\email{j.a.hoogeveen@uu.nl}

\affil*[1]{\orgdiv{Department of Information and Computing Science}, \orgname{Utrecht University}, \orgaddress{\street{Princetonplein 5}, \city{Utrecht}, \postcode{3584CC}, \country{Netherlands}}}

\affil[2]{\orgname{NS}, \orgaddress{\street{Laan van Puntenburg 100}, \city{Utrecht}, \postcode{3511ER}, \country{Netherlands}}}

\abstract{In many real-world scheduling problems, the processing times of tasks are subject to uncertainty. This makes it essential to design schedules that are robust against disruptions. We distinguish between two types of robustness: quality-robust schedules keep a stable objective value under disruptions, while solution-robust schedules keep a stable schedule structure. We study how to generate schedules for the Stochastic Parallel Machine Scheduling Problem (SPMSP) with precedence constraints and release dates that are both quality- and solution-robust. To this end, we compare several robustness measures and local search algorithms that use these measures. We start with an extensive simulation study where we determine the correlation of 13 existing and 6 new robustness measures with both quality and solution robustness. Next, we use the measures with high robustness correlation as the objective in our local search algorithm to find the best way to generate schedules containing both job-machine assignments as well as buffer assignments. We use a new dynamic objective weight within our local search algorithms, which actively changes the objective weights in each iteration. We found that the best way to generate robust schedules is to use robustness measures that model the processing times with normal distribution approximations and to combine job placement and buffer assignment simultaneously.}


\maketitle

\section{Introduction} \label{sec:intro}
Scheduling problems are among the most challenging classes of combinatorial optimization problems, as they typically involve allocating limited resources to tasks under complex temporal constraints. Their computational difficulty grows quickly with problem size, making exact solution methods impractical for many real-world applications. Over the past decades, local search has emerged as one of the most effective approaches for tackling such problems (see \citet{voss2012meta} for a survey of applications of local search). By iteratively improving candidate solutions through neighbourhood moves, local search methods have provided high-quality solutions within reasonable computational times across a wide variety of scheduling contexts, including production planning, timetabling, and transportation.

Many of the problems where local search has been most successful use deterministic data. However, disturbances are ubiquitous in real-world operations (\citet{xu2013robust}, \citet{allahverdi2010heuristics}). Due to these disturbances, deterministic models give an inaccurate representation of the objective value. Obtaining accurate objective values requires introducing stochastic factors into the model, for example, by replacing deterministic processing times with stochastic ones. However, for many problems, deriving a closed-form distribution of all stochastic variables is impossible, even when all distributions of individual processing times are known (see, for instance, \citet{passage2025new}). This intractability makes it difficult to apply local search methods, which typically rely on evaluating and comparing objective values directly. In this paper, we address this challenge in the context of the Stochastic Parallel Machine Scheduling Problem (SPMSP) with precedence constraints, release dates, and a global deadline. This problem has many real-world applications such as semiconductor fabrication, GPU scheduling, and maintenance scheduling. It presents an interconnection of constraints that makes the computational intractability especially acute.

Schedule performance under stochastic processing times is often estimated via simulation over many realizations, an approach adopted by \citet{van2013finding}. However, these simulations are computationally expensive, making them infeasible for local search algorithms, which require many iterations per second. This calls for a measure that is faster to evaluate while still indicating a schedule's robustness to disturbances — a role that surrogate robustness measures are designed to fill: given a complete schedule, they can be calculated efficiently and should reflect its robustness.

Many such measures exist in the literature, and most rely on a common mechanism: incorporating slack at the right points in the schedule to absorb disturbances before they propagate. Other measures, however, do not rely on slack directly; instead, they capture robustness through different structural or statistical properties of the schedule. There is no consensus on which of these measures, slack-based or otherwise, are best suited to producing robust schedules. To address this gap, we identify 13 existing robustness measures, for example from \citet{canon2009evaluation} and \citet{khemakhem2013efficient}, and introduce 6 new ones. We then conduct extensive simulation studies to compare these 19 measures and identify which are most strongly correlated with schedule robustness.

For these correlations, we discern two types of robustness. These are \emph{quality robustness}, which we define as the stability of our objective function under the effect of disturbances (i.e., whether the disturbances do not cause the makespan to exceed the deadline), and \emph{solution robustness}, which we define as the stability of the schedule itself under the effect of disturbances (i.e., whether all tasks can still start according to their original starting time). Although the literature focuses more on quality robustness, solution robustness is also vital. For instance, when jobs must be coupled to specific crew shifts, all jobs must start on time and can not be completed outside those shifts. This distinction also affects where slack should be inserted: buffering against deadline violations calls for different slack placement than buffering against shifts in individual start times.

After determining which measures have the highest correlation with schedule robustness, we also want to find the best way to generate robust schedules. An important aspect of generating robust schedules is deciding how to distribute slacks in a schedule. Multiple strategies for distributing slack in a schedule have been proposed in the literature; however, no conclusive comparison showing which method is best exists. Therefore, we will also compare different algorithms for distributing slack in a schedule. For these algorithms, we apply a dynamic solution weight that changes each iteration, balancing quality and solution robustness. We will also use a selection of our best quality robustness measures and solution robustness measures to definitively decide which are best for creating robust schedules. 

The rest of this paper is organised as follows. We begin with an overview of the existing literature on robustness measures in Section \ref{sec:lit}. Next, we provide a more in-depth explanation of our robustness measures in Section \ref{sec:theo-back}. In Section \ref{sec:Measure-correlation}, we present the setup and results of our first experiment, in which we determine the correlation between our robustness measures and the actual robustness of our schedules. In Section \ref{sec:buffer-insertion}, we use the best-performing measures from Section~\ref{sec:Measure-correlation} to compare local search strategies that distribute slack in different ways. Finally, we summarize our results and provide concluding remarks in Section~\ref{sec:conc}.

\subsection{Our Contribution}

In this paper, we investigate how to generate robust schedules for the SPMSP by evaluating and applying surrogate robustness measures within a local search algorithm. We study both quality robustness and solution robustness, examining how slack insertion strategies contribute to achieving each. Our paper presents four main contributions. 
\begin{enumerate}
    \item We compare 13 existing and 6 new robustness measures, evaluating their correlation with both quality robustness and solution robustness for the SPMSP.
    \item We embed these measures as objectives in a local search algorithm, comparing their performance for generating schedules that are both quality-robust and solution-robust.
    \item We introduce a dynamic solution weight function to balance both types of robustness throughout the local search procedure. To the best of our knowledge, this is the first time such a dynamic solution weight has been applied for this setting.
    \item We compare different strategies for assigning buffers using local search algorithms. We compare assigning buffers and jobs sequentially with assigning them in an integrated fashion. We evaluate which strategy produces the most quality- and solution-robust schedules.
\end{enumerate}

The overarching theme of these contributions is to compare different types of methodologies for generating schedules that are both quality-robust and solution-robust. Although many such methods exist, there is no consensus on which methods are best for creating quality- and solution-robust schedules.

\section{Literature Overview}\label{sec:lit}
We will first review the available literature for surrogate robustness measures and solution algorithms for the SPMSP. We will also include research on the Resource-Constrained Project Scheduling Problem (RCPSP), as it is a closely related problem that brings up interesting related work.

\subsection{Surrogate Robustness Measures}

\begin{figure}
    \centering
    \includegraphics[width = 0.7\textwidth]{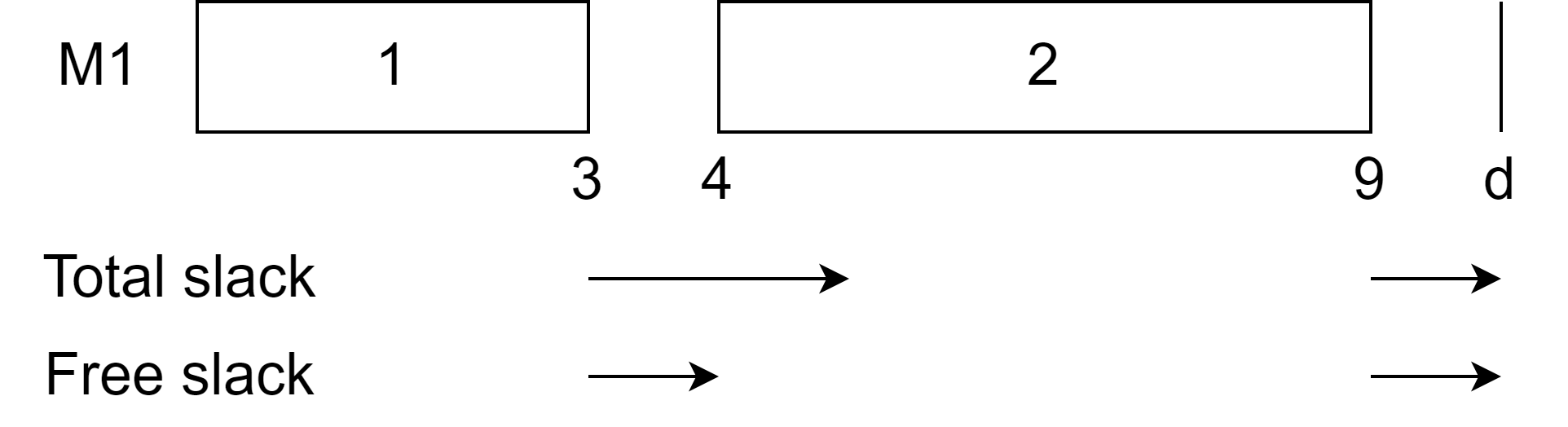}
    \caption{Illustration of the total slack and free slack for each job in a simple schedule. We see that for job 1 the total slack is larger than the free slack; job 2 can be delayed without any deadline violations. For job 2, the two kinds of slack are equal, as it is the final job in the schedule.}
    \label{fig:slack}
\end{figure}

Many existing measures in the literature are based on the concept of \emph{slack}. Here, a distinction is made between the \emph{total slack} of a job, which indicates how much a job can be delayed before the objective value (i.e., makespan) is changed, and the \emph{free slack} of a job, which indicates how much a job can be delayed before other jobs are impacted; see also Figure \ref{fig:slack}. Free slack is also commonly referred to as a \emph{buffer}. One of the first contributions is by \citet{jorge1994robustness}, who investigate robustness for the Job Shop Scheduling Problem. They propose the average total slack of all jobs as a robustness measure. They perform a simulation study which shows that the average total slack of a schedule is a good indicator of expected schedule delay. \citet{hazir2010robust} investigate robustness measures for the Discrete Time/Cost Trade-Off Problem using Monte Carlo simulation. They find that the measure that is most correlated with their performance metrics is the project buffer size, which is similar to the minimum total slack. \citet{al2005bi} propose using the sum of free slacks as a robustness measure for the RCPSP, rather than the minimum total slack. This, however, is critiqued by \citet{kobylanski2007note}, who show that using the total free slack can lead to schedules that are far from robust. Instead, they propose using the schedule's minimum free slack as a robustness measure. Additionally, they argue that for practical applications, it is better to fix the makespan and then determine a robust schedule for that makespan. \citet{chtourou2008two} also use the free slack to increase the robustness of the schedule for the RCPSP. They do this by proposing a two-stage priority-rule-based approach. They first solve the problem of minimizing the makespan, and then use free-slack-based measures to maximize robustness while keeping the makespan at that minimum.

Multiple studies propose robustness measures based on probability distributions. \citet{canon2009evaluation} investigate multiple such measures. They introduce the standard deviation of the makespan distribution and the makespan 0.99-quantile as robustness measures. They also introduce the slack mean, which is the expected value of the sum of the free slacks. They compare these measures with some already existing measures. The first of these is the makespan differential entropy, introduced by \citet{boloni2002robust}. This measure aims to quantify the variance of the distribution, which they argue gives a good indication of the robustness. Secondly, they consider a measure introduced by \citet{shestak2006stochastic}, which quantifies the probability that the makespan is between two bounds. These can be absolute or relative bounds. The last measure included is the lateness likelihood introduced by \citet{shi2006robust}, which is defined as the probability that the makespan exceeds a chosen target. \citet{canon2009evaluation} compare these measures in a large simulation study in which they assume they know all processing time distributions. They show that almost all measures are equivalent. Therefore, they argue that the simplest measure, the makespan standard deviation, will be sufficient in most cases. They also find that the sum of total slacks has at most a weak correlation with all other quality robustness measures. Another way to approximate the makespan using probability distributions has been presented by \citet{van2013finding}. They use discrete-event simulation to approximate the makespan, and show that it outperforms some classic methods. Lastly, \citet{passage2025new} recursively estimate the expected value and variance of each start and completion time, approximating the schedule makespan. A key element in this approach is determining the expected value and variance of the maximum of two normal distributions. For these estimations, they use the methods of \citet{nadarajah2008exact}. They found that it performs very well while being relatively quick to compute.

Given this variety of measures, researchers have compared them to discover which measures generally perform best. Firstly, there is the study by \citet{canon2009evaluation} that we mentioned in the previous paragraph. There is also the study by \citet{khemakhem2013efficient} which compares the performance of several measures proposed by \citet{al2005bi}, \citet{kobylanski2007note}, \citet{lambrechts2008tabu} and \citet{chtourou2008two} for the RCPSP. Besides these known measures, they introduce some new ones. These new measures are based on \emph{slack sufficiency}, which counts for every job whether the available slack is sufficient to absorb its expected delay. They extend this by also including the slack sufficiency of preceding jobs, and find that their newly proposed measures outperform existing ones. They further find that the best-performing measures account for the possible delay of both a job itself and its predecessors. Another relevant study is that of \citet{van2018measure}, who investigate robustness measures for shunting schedules, which can be modelled as RCPSPs. They compare existing measures from \citet{jorge1994robustness}, \citet{al2005bi}, \citet{kobylanski2007note}, \citet{khemakhem2013efficient}, \citet{wilson2014flexibility}, and \citet{passage2025new}. Additionally, they propose one new measure adapted from \cite{wilson2014flexibility} and one based on the slack of a path of jobs in the schedule; by assuming that the jobs in that path are normally distributed, they can estimate the makespan of that path. In the end, they show that robustness measures based on approximation with a normal distribution, such as the newly proposed one and the one proposed by \citet{passage2025new}, strongly correlate with the fraction of simulations that do not meet the deadline, and the average schedule lateness. For this study, they work with earliest start schedules, meaning that a job can start as soon as its predecessors are finished.

Most papers mentioned in this section base the performance of their robustness measures on the quality robustness of the resulting schedules. Of the ones we mentioned so far, only \citet{kobylanski2007note} investigate solution robustness. They find that maximizing the minimum ratio of free slack to activity duration increases both quality and solution robustness. \citet{van2005use} argue that an ideal schedule should combine both quality robustness and solution robustness. For a closer look at solution robustness itself, we refer the reader to \citet{ke2015uncertain}, who use the expected deviation between the realized and planned start time to create robust schedules for the RCPSP. Another relevant contribution is that of \citet{fu2015robust}, who apply a robust optimization framework to improve solution robustness for the Generalized Assignment Problem. However, there has not been a good comparative study that investigates which measures work well for solution robustness.

Despite many studies investigating robustness measures and their effectiveness, there is no clear consensus on which type of measure is useful for which situation. While the goal is often to ensure a schedule is completed before the deadline, it can also be very important that all tasks in the schedule start at their expected starting time. Therefore, additional research is required to discover which type of measures work well for quality robustness and which for solution robustness. In this paper, we introduce 6 new measures, and also use measures from \citet{jorge1994robustness}, \citet{al2005bi}, \citet{kobylanski2007note}, \citet{chtourou2008two}, \citet{khemakhem2013efficient}, \citet{wilson2014flexibility}, \citet{van2018measure}, and \citet{passage2025new}.

\subsection{Algorithms for creating robust schedules}

The literature contains a variety of methods for creating robust schedules for the SPMSP and RCPSP. An important part of robust scheduling is the method used for buffer assignment. We start by reviewing some methods that use sequential phases for job and buffer assignment. First, there is \citet{lambrechts2008proactive}, who developed an algorithm for creating robust schedules for the RCPSP. They use sequential phases for adding resource buffers and time buffers. This method has been used and expanded by \citet{lambrechts2011time} and \citet{ma2019computational}, who found it effective for creating robust baseline schedules for the RCPSP. \citet{bruni2011heuristic} used a different heuristic method; they first schedule all jobs in a rolling horizon fashion, and afterwards similarly insert buffers. The size of the buffers is determined using probability chance constraints. This approach allows for choosing the desired risk level at different points in the schedule.  Next, \citet{liang2020integration} use a branch-and-bound method to find a schedule with minimum deterministic makespan. Afterwards, buffers are inserted using simulated annealing. They use a bi-objective model to optimize both the quality robustness and solution robustness. The authors found that although their method is able to find robust schedules within low computation time, it is unclear how much improvement there would be from using a method that integrates job and buffer allocation.

Several papers used an integrated method for allocating jobs and buffers. Firstly, \citet{skutella2005stochastic} define a list scheduling policy for the SPMSP with precedence constraints. Although they do not test the performance of this algorithm, they derive theoretical performance guarantees, focusing on the quality robustness of the schedule. \citet{lambrechts2008tabu} use a tabu search algorithm to find a schedule for the RCPSP with maximum solution robustness. They wanted the schedule to be robust against uncertainties in the maximum resource usage at every time. They compared their method to several other methods that have separate phases for job and buffer assignments and found that their method outperformed most of these. \citet{ma2022genetic} also looked mainly at solution robustness for the RCPSP, but instead used a genetic algorithm for this. This genetic algorithm also integrated the job and buffer assignment. They compared their method with \citet{lambrechts2008tabu} and found that their genetic algorithm produced more robust schedules more quickly. We also found two papers that focused on the SPMSP without precedence constraints. \citet{liu2021parallel} used a stochastic programming approach to account for stochastic release dates and processing times. They found that their method, which uses scenario reduction, outperformed sample approximation, but they still hit the computation time limit of one hour for randomly generated instances with 9 or more jobs. \citet{xiao2026branch} used a branch-cut-and-price algorithm to solve the SPMSP with release dates, but did not include precedence constraints. They found that their method produced high-quality solutions for instances with up to 20 machines and 120 jobs.

When we compare the objectives used in the literature discussed above, we see that most work does not consider optimizing both quality and solution robustness together. Most work focuses either on one of the two \citep{bruni2011heuristic,liu2021parallel, skutella2005stochastic, xiao2026branch} or fixes quality robustness with a deadline constraint and then optimizes solution robustness \citep{lambrechts2008proactive, lambrechts2008tabu, lambrechts2011time, ma2019computational, ma2022genetic}. Only one paper in this body of work explicitly investigates a trade-off between the two \citep{liang2020integration}.

Although we see that there is already some research comparing different types of algorithms for creating buffered schedules, there is currently no existing work that considers an integrated approach that optimizes for both quality and solution robustness. Moreover, there is also no existing literature that includes multiple measures for both quality and solution robustness. We aim to close this gap.

\section{Robustness Measures}\label{sec:theo-back}
This section describes our problem setting, the Stochastic Parallel Machine Scheduling Problem, and defines the robustness measures we investigate. In the SPMSP, we have $m$ identical machines that have to carry out $n$ jobs. Each machine can process at most one job at the same time. The processing time of each job $j$ ($j = 1,\dots, n$) is a stochastic variable denoted by $D_j$. The average of this distribution will be denoted by $p_j$. Each job $j$ also has a release date $r_j$ before which it cannot start. We consider a global deadline constraint for this problem, meaning that each job must be completed before a global deadline $d$. Additionally, precedence constraints between jobs give us a partial ordering that determines when jobs can be executed. When we have a precedence constraint $i\prec j$, job $i$ must finish before job $j$ can start. The goal is then to find a schedule that is feasible under all these constraints and is as robust as possible. 

In this context, a schedule describes the (planned) starting time of each job and the machine on which the job is processed. For a given schedule $\sigma$, we denote the starting time of job $j$ by $s_j^{\sigma}$. For each job $j$ ($j = 1,\dots,n$), we also compute a latest starting time $ls_j^{\sigma}$: the latest time the job can start such that the schedule makespan does not exceed the deadline, while respecting all precedence constraints. Because processing times are stochastic, this latest starting time cannot be guaranteed with certainty; we therefore compute $ls_j^{\sigma}$ under the assumption that each job takes its expected processing time $p_j$. When simulating the schedule, however, we use the full processing time distributions rather than these expected values. In addition, during simulation, a job is never allowed to start before its planned starting time $s_j^{\sigma}$, even if its predecessors happen to finish earlier than planned. This execution policy will be used throughout the whole paper.

\subsection{Slack based measures}

As mentioned in the introduction, a common way to make schedules robust against disruptions is by including slack. For our robustness measures, we differentiate between two types of slack. Firstly, there is the \textit{total slack} of a job $j$ given by 
\begin{equation}
ts_j^\sigma = ls_j^\sigma - s_j^\sigma.
\end{equation}
Secondly, there is also the free slack $fs_j^\sigma$. This is the amount of time that job $j$ can be delayed without impacting any of its successor jobs. This is calculated as:
\begin{equation}
    fs_j^\sigma \coloneq \min_{i\in\operatorname{dsucc_\sigma(j)}}\{s_i^\sigma - s_j^\sigma - p_j\}.
\end{equation}
Here $\operatorname{dsucc}_\sigma(j)$ denotes all of the direct successor jobs of job $j$. Note that when we are talking about (direct) predecessors or successors, we will always be talking about both job predecessors and successors as well as machine predecessors and successors unless specified otherwise. When we have two jobs $i, j$, where $i$ is a predecessor of $j$, we will write $i\prec_\sigma j$; we include the $\sigma$ subscript to indicate that this is specifically in schedule $\sigma$ and hence also includes machine predecessors. Furthermore, when we talk about a \emph{direct} predecessor of a job $j$, we mean any job $i$ with $i\prec_\sigma j$ such that there is no job $k$ with $i\prec_\sigma k\prec_\sigma j$. A complete overview of all notation used is found in Table \ref{tab:notation}.

\begin{table}[htbp]
    \centering
    \caption{Overview of all notation used in this paper.}
    \begin{tabular}{lp{11cm}}
    \toprule
    Variable & Definition \\ \midrule
    $\sigma$ & A schedule.\\
    $n$ & The number of jobs.\\
    $m$ & The number of machines.\\
    $j$ & A job.\\
    $p_j$ & Expected processing time of job $j$.\\
    $D_j$ & Stochastic variable indicating the processing time distribution for job $j$.\\
    $\sigma_j$ & Standard deviation of $D_j$.\\
    $r_j$ & Release date of job $j$.\\
    $d$ & Deadline of the schedule.\\
    $i\prec_\sigma j$ & Job $i$ is a predecessor of job $j$ in schedule $\sigma$, meaning job $i$ must be completed before job $j$.\\
    $s_j^\sigma$ & Planned starting time of job $j$ in schedule $\sigma$.\\
    $ls_j^\sigma$ & Latest possible starting time for job $j$ in schedule $\sigma$ based on expected processing time $p_j$.\\
    $ts_j^\sigma$ & Total slack of job $j$ in schedule $\sigma$.\\
    $fs_j^\sigma$ & Free slack of job $j$ in schedule $\sigma$.\\
    $ndp_j^\sigma$ & Number of direct predecessors of job $j$ in schedule $\sigma$.\\
    $nds_j^\sigma$ & Number of direct successors of job $j$ in schedule $\sigma$.\\
    $\operatorname{prec}_j^\sigma$ & Set of all predecessors of job $j$ in schedule $\sigma$.\\
    $\operatorname{dprec}_j^\sigma$ & Set of the direct predecessors of job $j$ in schedule $\sigma$.\\
    $\operatorname{succ}_j^\sigma$ & Set of all successors of job $j$ in schedule $\sigma$.\\
    $\operatorname{dsucc}_j^\sigma$ & Set of the direct successors of job $j$ in schedule $\sigma$.\\
    $X_j^\sigma$ & Stochastic variable indicating the starting time distribution for job $j$ in schedule $\sigma$.\\
    $Y_j^\sigma$ & Stochastic variable indicating the completion time distribution for job $j$ in schedule $\sigma$.\\
    $ESD_j^\sigma$ & Estimated starting delay of job $j$ in schedule $\sigma$.\\\bottomrule
    \end{tabular}
    \label{tab:notation}
\end{table}

A comparison between the two types of slack is shown in Figure \ref{fig:slack}. In the figure, we see that both job 1 and job 2 have a small buffer between their completion times and either the starting time of the next job (for job 1) or the deadline (for job 2). The size of this buffer equals the free slack. We also see that job 1 has a higher total slack, because we can still delay job 2 without violating the deadline. So free slack is like a buffer before the next job begins, while total slack is the total extra buffer time between a job's end and its deadline (taking all successor jobs into account).

Our first two robustness measures are to sum either all total slacks or all free slacks in the schedule. These measures have been proposed earlier by \citet{jorge1994robustness} and \citet{al2005bi}.
\begin{align}
    RM_1(\sigma) &= \sum_{j\in[n]}ts_j^\sigma,\\
    RM_2(\sigma) &= \sum_{j\in[n]}fs_j^\sigma,
\end{align}
where $j\in[n]$ is shorthand notation for $j\in\{1,\dots,n\}$. Next, we have two measures proposed by \citet{kobylanski2007note}. The idea behind these two measures is that a schedule will be more robust if the slack is more evenly spread between the different jobs.
\begin{align}
    RM_3(\sigma) &= \min_{j\in[n]}\{ts_j^\sigma\},\\
    RM_4(\sigma) &= \min_{j\in[n]}\{fs_j^\sigma/p_j\}.
\end{align}
The fifth measure is one proposed in \citet{chtourou2008two}. It also counts the sum of the free slacks, but for each summand, it takes the minimum between the free slack and some multiple of the processing time:
\begin{align}
    RM_5(\sigma) = \sum_{j\in[n]} \min\{fs_j^\sigma, \lambda_j p_j\}.
\end{align}
This measure uses a parameter $0 \leq \lambda_j\leq 1$ for each job $j\in [n]$. \citet{chtourou2008two} suggest choosing $\lambda_j$ such that it equals the expected fractional deviation from the average processing time. So if we have some job $j$ that has distribution $D_j$, then the recommended value for $\lambda_j$ is given by $\lambda_j = \frac{1}{p_j}\EE(\abs{D_j - p_j})$. The idea is that summing the minimum of the expected delay and the free slack, any slack that might wrongfully inflate the total free slack is ignored. Any future use of the variable $\lambda_j$ in this section uses this same idea. Alternatively, they also propose to count how often the free slack is non-zero:
\begin{align}
    RM_6(\sigma) &= \sum_{j\in[n]}\alpha_j, & \text{where }\alpha_j = \begin{cases}
        1 & fs_j^\sigma > 0,\\0 & \text{otherwise}.
    \end{cases}
\end{align}
Next, we consider several measures where, like earlier, we take the sum of all free slacks. This time, however, we will be weighing the terms in the sum by different values. For these weights, \citet{chtourou2008two} proposed using the expected duration ($p_j$) of the job, the number of direct successors ($nds_j^\sigma$) of a job, and the product of the expected duration and the number of direct successors. As a new measure, we propose to use the number of direct predecessors ($ndp_j^\sigma$) of a job. These give us the following four measures:
\begin{align}
    RM_7(\sigma) &= \sum_{j\in[n]}fs_j^\sigma\times p_j,\\
    RM_8(\sigma) &= \sum_{j\in[n]}fs_j^\sigma\times ndp_j^\sigma,\\
    RM_9(\sigma) &= \sum_{j\in[n]}fs_j^\sigma\times nds_j^\sigma,\\
    RM_{10}(\sigma) &= \sum_{j\in[n]}fs_j^\sigma\times nds_j^\sigma\times p_j.
\end{align}
For the next two measures, we consider how well each job can absorb the delay of its predecessors. \citet{khemakhem2013efficient} propose, for each job $j$, to count the number of predecessors whose expected delay is smaller than the free slack of $j$ --- including its own expected delay. Here, ``predecessors'' refers not just to direct predecessors, but to all jobs that must finish before $j$ can start. This count reflects how well a job can recover from delay accumulated earlier in the schedule: jobs with too little slack relative to their predecessors' delays are the most critical for robustness, and jobs with many successors are implicitly weighted more heavily, since they appear in the count for each of those successors. We also propose counting the reverse (i.e., the number of predecessors with \emph{insufficient} slack), in which case a lower value indicates a more robust schedule. Below, we use $\operatorname{prec}_j^\sigma$ to denote the set of all predecessors of job $j$ in schedule $\sigma$.
\begin{align}
    RM_{11}(\sigma) &= \sum_{j\in[n]}\#\{i\mid i\in \operatorname{prec}_j^\sigma\cup\{j\}, fs_j^\sigma\geq \lambda_i p_i\}.\\
    RM_{12}(\sigma) &= \sum_{j\in[n]}\#\{i\mid i\in \operatorname{prec}_j^\sigma\cup\{j\}, fs_j^\sigma< \lambda_i p_i\}.
\end{align}
The reason we also introduced the reverse is that schedules with more precedence constraints might give a higher value for $RM_{11}$ even though those extra precedence constraints might have made the schedule less robust overall. This problem does not exist for $RM_{12}$.
\subsection{Interval schedules}
\citet{wilson2014flexibility} propose a measure based on linear programming, where each job $j$ is given an interval $[e_i, l_i]$ representing the possible start times of the job. The objective is then to maximize the sum of all of the intervals. This linear program is given below:
    \begin{maxi}
        {}{\sum_{j\in[n]}\qty(l_j - e_j)}{}{RM_{13}(\sigma) = }
        \addConstraint{s_j^\sigma}{\leq e_j\leq l_j}{\quad\forall j}
        \addConstraint{l_j + p_j}{\leq e_i}{\quad\forall j\in\operatorname{prec}_i^{\sigma}}
        \addConstraint{l_j,e_j}{\leq d - p_j}{\quad\forall j\in[n]}
    \end{maxi}
In the same vein, \citet{van2018measure} used the same linear program with a different objective, aiming to spread the intervals more evenly across different jobs.
    \begin{maxi}
        {}{\min_{j\in[n]}\{l_j - e_j\}}{}{RM_{14}(\sigma) = }
        \addConstraint{s_j^\sigma}{\leq e_j\leq l_j}{\quad\forall j}
        \addConstraint{l_j + p_j}{\leq e_i}{\quad\forall j\in\operatorname{prec}
        _i^{\sigma}}\label{eq:RM14}
        \addConstraint{l_j,e_j}{\leq d - p_j}{\quad\forall j\in[n]}
    \end{maxi}

\subsection{Normal approximation}
Next, we will consider three measures that use an approximation for the maximum of two normal distributions. It is known that the maximum and minimum of two normal distributions are, in general, not normally distributed (\cite{nadarajah2008exact}). However, \citet{passage2025new} propose a method to make a very good approximation of this distribution using a normal distribution. This allows us to approximate the complete makespan distribution by continuously applying this approximation for the maximum of two normal distributions. We let $X_j^\sigma$ denote the distribution for the starting time of job $j$ in schedule $\sigma$, and we let $Y_j^\sigma$ denote the distribution for the completion time of job $j$ in schedule $\sigma$. Then we can calculate $X_j^\sigma$ and $Y_{j}^\sigma$ as 
\begin{align}
    X_j^\sigma &= \max\Big\{\max_{i\in\operatorname{dprec}_j^\sigma}\{Y_i^\sigma\}, s_{j}^\sigma\Big\}\\
    Y_{j}^\sigma &= X_j^\sigma + D_{j},
\end{align}
where $D_j$ denotes the random variable for the duration of job $j$. In words, to calculate $X_j^\sigma$, we take the maximum completion time of all direct predecessors of job $j$, including the machine predecessor. This is done by sequentially taking the maximum between the current maximum and the next distribution to add. After this, we take the maximum of that time with the planned starting time of job $j$, from which we obtain the starting time distribution of job $j$. Then, we calculate $Y_j^\sigma$ by adding the duration distribution for job $j$. This assumes that we have estimates of the distributions for all predecessors. We calculate these using a topological ordering of all of the jobs. If a job has no predecessors, then of course the distributions are easy to calculate. With these definitions, we find three more robustness measures. The first of which can also be found in \citet{passage2025new}, while the second and third are two of our newly proposed measures.
\begin{align}
    RM_{15}(\sigma) &= \PP(\max_{j\in[n]}Y_j^\sigma\leq d),\\
    RM_{16}(\sigma) &= \sum_{j\in[n]}\PP\qty(\max_{i\in\operatorname{dprec}_j^\sigma}\{Y_i^\sigma\}\leq s_j^\sigma),\\
    RM_{17}(\sigma) &= \sum_{j\in[n]}\mathbb{E}(Y_j^\sigma) - s_j^\sigma
\end{align}
Here $d$ denotes the deadline for the schedule. In words, the first of these three looks at the probability that the schedule finishes on time. The second measure calculates the sum over all jobs, where for each job we sum the probability that its predecessors finish before its starting time (and hence the job can start on time). Thus, if this value is high, we expect most jobs to start on time, indicating a high solution robustness. For the last measure, we sum the difference between the expected and planned starting time for each job. Thus, a higher value would indicate a less robust schedule.

\subsection{Predecessor slack}
Lastly, we have two new measures that are based on the slack of predecessor jobs. The idea here is that the probability that a job starts on time completely depends on the slack of its predecessors. The first of the two measures is defined as:
\begin{align}
    RM_{18}(\sigma) = \sum_{j\in[n]} \beta_j,
\end{align}
where
\begin{align}
    \beta_j\coloneq \begin{cases} \frac{\#\{i\mid fs_i^\sigma\geq \lambda_i p_i, i\in \operatorname{dprec}_j^\sigma\}}{\#\operatorname{dprec}_j^\sigma} & \text{if }\#\operatorname{dprec}_j^\sigma > 0,\\ 1 & \text{otherwise}.\end{cases}
\end{align}
Here, $\operatorname{dprec}_j^\sigma$ denotes the set of direct predecessors of job $j$ in $\sigma$. The idea behind this measure is that the quantity $\beta_j$ indicates how well the predecessors of job $j$ can handle delays. This thus indicates the likelihood of job $j$ starting on time, meaning a higher value should correspond to a more robust schedule. The main difference between this new measure and the previously defined $RM_{11}$ and $RM_{12}$ is that we consider the free slack of the jobs directly preceding the current job instead of the free slack of the job itself. For the last measure, we will use an estimate of the starting delay (ESD) of a job, which is defined as:
\begin{align}
    ESD_j^\sigma \coloneq \begin{cases}
        0 & \operatorname{prec}_j^\sigma = \emptyset\\
        \max_{i\prec_\sigma j}\{\max\{\lambda_i p_i + ESD_i^\sigma - fs_i^\sigma, 0\}\} & \text{otherwise}
    \end{cases}.
\end{align}
For the robustness measure, we then sum over all of the starting delays.
\begin{align}
    RM_{19}(\sigma) = \sum_{j\in[n]}ESD_j^\sigma.
\end{align}
This value can be calculated using dynamic programming; the schedule is topologically sorted and then traversed in that order. Of course, for this measure, a higher value means we expect more delay and hence a less robust schedule. Note that for both $RM_{18}$ and $RM_{19}$ we again use a constant $\lambda_i p_i$. Just as before, we can think of this constant as the expected increase in processing time.

When we look at all of the measures that we have described in this section, it might seem like some of these are quite difficult to compute. However, most of the measures in this section require us to traverse the schedule at most 1 time, making them significantly cheaper to compute than simulation and very suitable for local search. We have also created an overview of all our robustness measures in Table \ref{tab:rob-measures}.

\begin{table}[htbp]
    \centering
    \caption{Overview of all robustness measures including formulae and short descriptions.}
    \begin{tabular}{lp{0.4\linewidth}p{0.35\linewidth}}
        \toprule
        \parbox[c]{2cm}{\raggedright Robustness Measure} & Formula & Description\\
        \midrule
        $RM_1$ & $\sum_{j\in[n]}ts_j^\sigma$ & Sum of total slacks.\\
        $RM_2$ & $\sum_{j\in[n]}fs_j^\sigma$ & Sum of free slacks.\\
        $RM_3$ & $\min_{j\in[n]}\{ts_j^\sigma\}$ & Minimum total slack.\\
        $RM_4$ & $\min_{j\in[n]}\{fs_j^\sigma/p_j\}$ & Minimum free slack weighted by job duration.\\
        $RM_5$ & $\sum_{j\in[n]} \min\{fs_j^\sigma, \lambda_j p_j\}$ &Sum of minima of free slack and expected delay.\\
        $RM_6$ & $\sum_{j\in[n]}\alpha_j$, $\alpha_j$ Indicates whether job $j$ has non-zero free slack & Number of jobs with non-zero free slack.\\
        $RM_7$ & $\sum_{j\in[n]}fs_j^\sigma\times p_j$ & Sum of free slacks weighted with job duration.\\
        $RM_8$ & $\sum_{j\in[n]}fs_j^\sigma\times ndp_j^\sigma$ & Sum of free slacks weighted with number of predecessors.\\
        $RM_9$ & $\sum_{j\in[n]}fs_j^\sigma\times nds_j^\sigma$ & Sum of free slacks weighted with number of successors.\\
        $RM_{10}$ & $\sum_{j\in[n]}fs_j^\sigma\times nds_j^\sigma\times p_j$ & Sum of free slacks weighted with number of successors and job duration.\\
        $RM_{11}$ & $\sum_{j\in[n]}\#\{i\mid i\in \operatorname{prec}_j^\sigma\cup\{j\}, fs_j^\sigma\geq \lambda_i p_i\}$ & Sum of the number of predecessors that have sufficient slack for all jobs.\\
        $RM_{12}$ & $\sum_{j\in[n]}\#\{i\mid i\in \operatorname{prec}_j^\sigma\cup\{j\}, fs_j^\sigma< \lambda_i p_i\}$ & Sum of the number of predecessors that have insufficient slack for all jobs.\\
        $RM_{13}$ & $\max \sum_j (l_j - e_j)$, $[e_j;l_j]$ execution interval for job $j$& Optimal value of LP that maximizes the sum of execution intervals.\\
        $RM_{14}$ & $\max \min_j (l_j - e_j)$, $[e_j;l_j]$ execution interval for job $j$ & Optimal value of LP that maximizes the minimum execution intervals.\\
        $RM_{15}$ & $\PP(\max_{j\in[n]}Y_j^\sigma\leq d)$ & Approximated probability of finishing on time.\\
        $RM_{16}$ & $\sum_{j\in[n]}\PP\qty(\max_{i\in\operatorname{dprec}_j^\sigma}\{Y_i^\sigma\}\leq s_j^\sigma)$ & Sum of approximated probabilities of starting on time.\\
        $RM_{17}$ & $\sum_{j\in[n]}\mathbb{E}(Y_j^\sigma) - s_j^\sigma$ & Expected starting delay based on normal approximation.\\
        $RM_{18}$ & $\sum_{j\in[n]} \beta_j$ & Sum of fractions of predecessors with sufficient free slack\\
        $RM_{19}$ & $\sum_{j\in[n]}ESD_j^\sigma$ & Sum of expected starting delays\\
        \bottomrule
    \end{tabular}
    \label{tab:rob-measures}
\end{table}

\section{Evaluation of the Robustness Measures for the SPMSP}\label{sec:Measure-correlation}
Before investigating the best method for generating robust schedules, we first determine which robustness measures correlate most strongly with actual schedule robustness and are thereby good candidates for quality or solution robustness measures. To this end, we generate a diverse set of schedules for multiple SPMSP instances and simulate their execution. We use a diverse set of instances with a variable number of jobs (30 to 100), machines (4 to 12), and precedence constraints (15 to 250). Appendix~\ref{sec:app-instance} describes how we obtain these instances in more detail. To ensure that we generate a diverse set of schedules for these instances, we first generate 10 baseline schedules without buffers for every instance. Then, we use a buffer insertion method to generate 97 different buffered schedules per baseline schedule. Thus, we end up with 970 schedules per instance for a total of 24 instances. The schedule generation is described in more detail in Appendix~\ref{sec:app-scehdules}.

\begin{figure}[htbp]
    \centering
    \includegraphics[width=\linewidth]{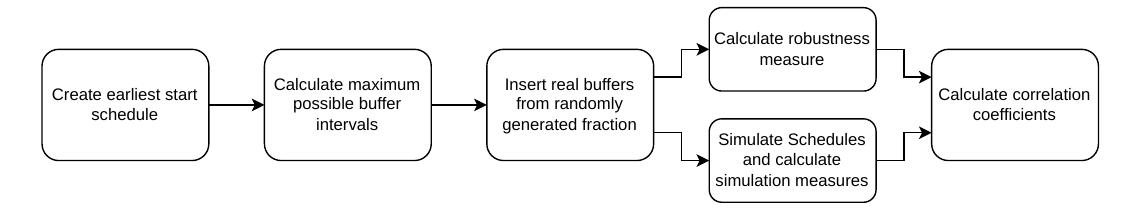}
    \caption{Schematic overview of the steps required to calculate the correlation coefficients for the robustness measures.}
    \label{fig:correlation-flow-diagram}
\end{figure}

For all these schedules, we can determine the values of the different measures. However, to determine which measures best indicate schedule robustness, we also need to simulate the schedules. For each schedule, we perform 1000 simulations. These simulations follow a strict execution policy, where jobs cannot be started before their planned starting time. Then, to assess the quality robustness of these schedules, we measure the average realized makespan and the average fraction of simulations completed within the deadline. For solution robustness we measure the average fraction of on-time jobs and the average total delay. An overview of these simulation measures can be found in Table~\ref{tab:simulation_measures}, and an overview of this whole process can be found in Figure~\ref{fig:correlation-flow-diagram}.

For the processing times of the jobs, we will use two continuous distributions: the normal distribution and the log-normal distribution. The mean of the distribution will always be the initial deterministic processing time, and we will set the standard deviation as some fraction of the mean. Specifically, for the normal distribution, we consider the fractions 0.25 and 0.4, while for the log-normal distribution we will consider 0.25. We refer to these distributions as $N_{25}$, $N_{40}$, and $LN_{25}$ respectively. Recall that we also need to choose values for the $\lambda_j$ parameters for some of our robustness measures. Although \citet{chtourou2008two} suggested using the expected increase in processing time, \citet{van2013finding} found that taking the 70th percentile often gave the best result; thus we will use this as well. The corresponding value of $\lambda_j$ can then be calculated from the inverse of the cumulative distribution function.

Now that we have calculated all robustness measure and simulation measure values, we use these to calculate the Spearman correlation coefficient $\rho$ for every pair of robustness measure and simulation measure. Figures~\ref{fig:corr_avgmakespan_initial}, \ref{fig:corr_fraction_initial}, \ref{fig:corr_onTime_initial}, and \ref{fig:corr_totaldelay_initial} report the spread of the absolute values of these coefficients, per robustness measure and per simulation measure. A value close to $1$ indicates good correlation between the robustness measure and the simulation measure, while a value close to $0$ indicates poor correlation.

\begin{table}[htbp]
    \centering
    \caption{Overview of our simulation measures for assessing the robustness of a schedule, and the type of robustness they assess for the SPMSP. To calculate these measures, we simulate a schedule 1000 times.}
    \begin{tabular}{cc}
    \toprule
         Quality Robustness & Solution Robustness \\
         \midrule
         Average makespan & Average fraction of on-time jobs\\
         Average fraction within deadline & Average total delay\\
         \bottomrule
    \end{tabular}
    \label{tab:simulation_measures}
\end{table}

\subsection{Results quality robustness}\label{sec:first-exp-results-quality}

\begin{figure}[htbp]
    \centering
    \includegraphics[width = \textwidth]{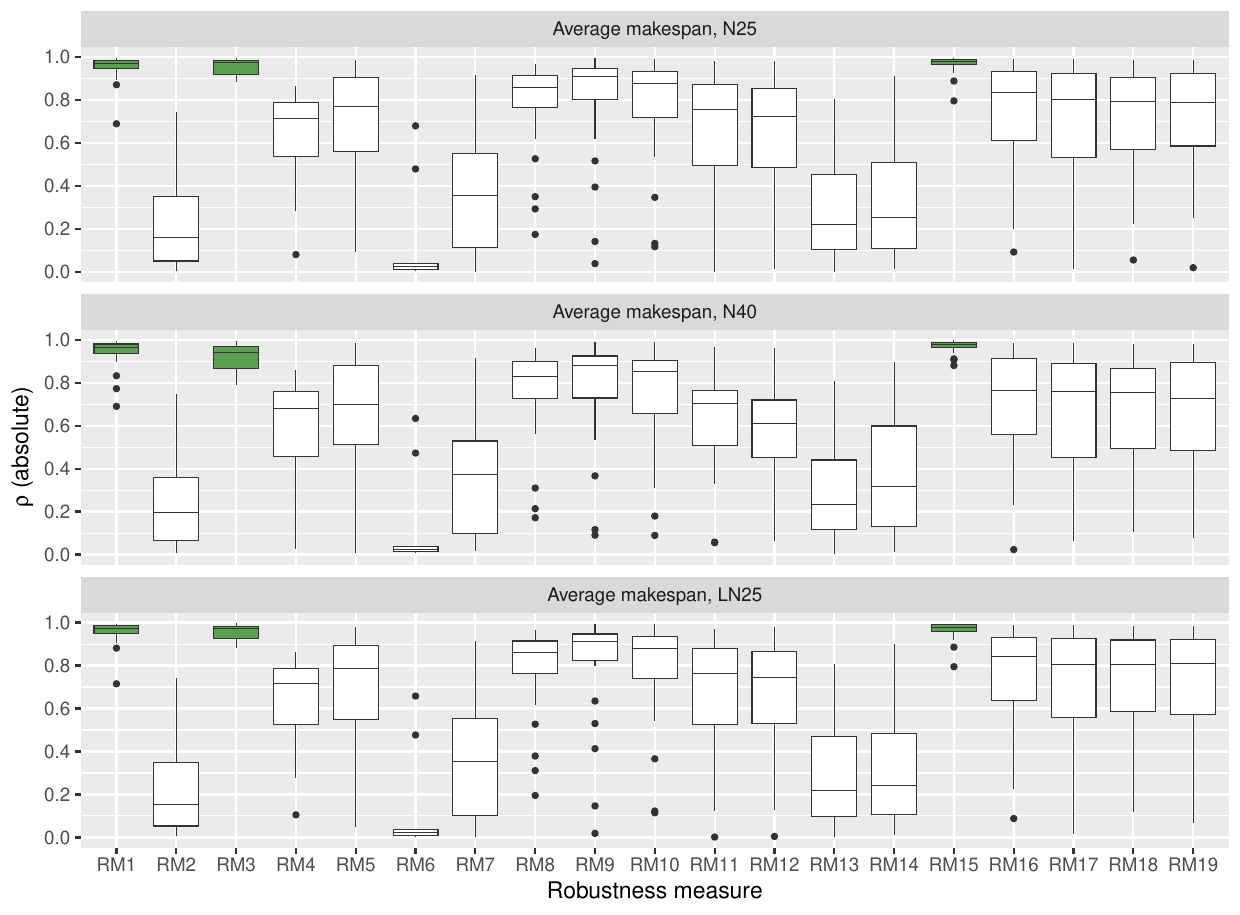}
    \caption{Spearman's coefficient $\rho$ (absolute value) for the average makespan. Green highlight indicates that the mean is greater than $0.9$.}
    \label{fig:corr_avgmakespan_initial}
\end{figure}

\begin{figure}[htbp]
    \centering
    \includegraphics[width = \textwidth]{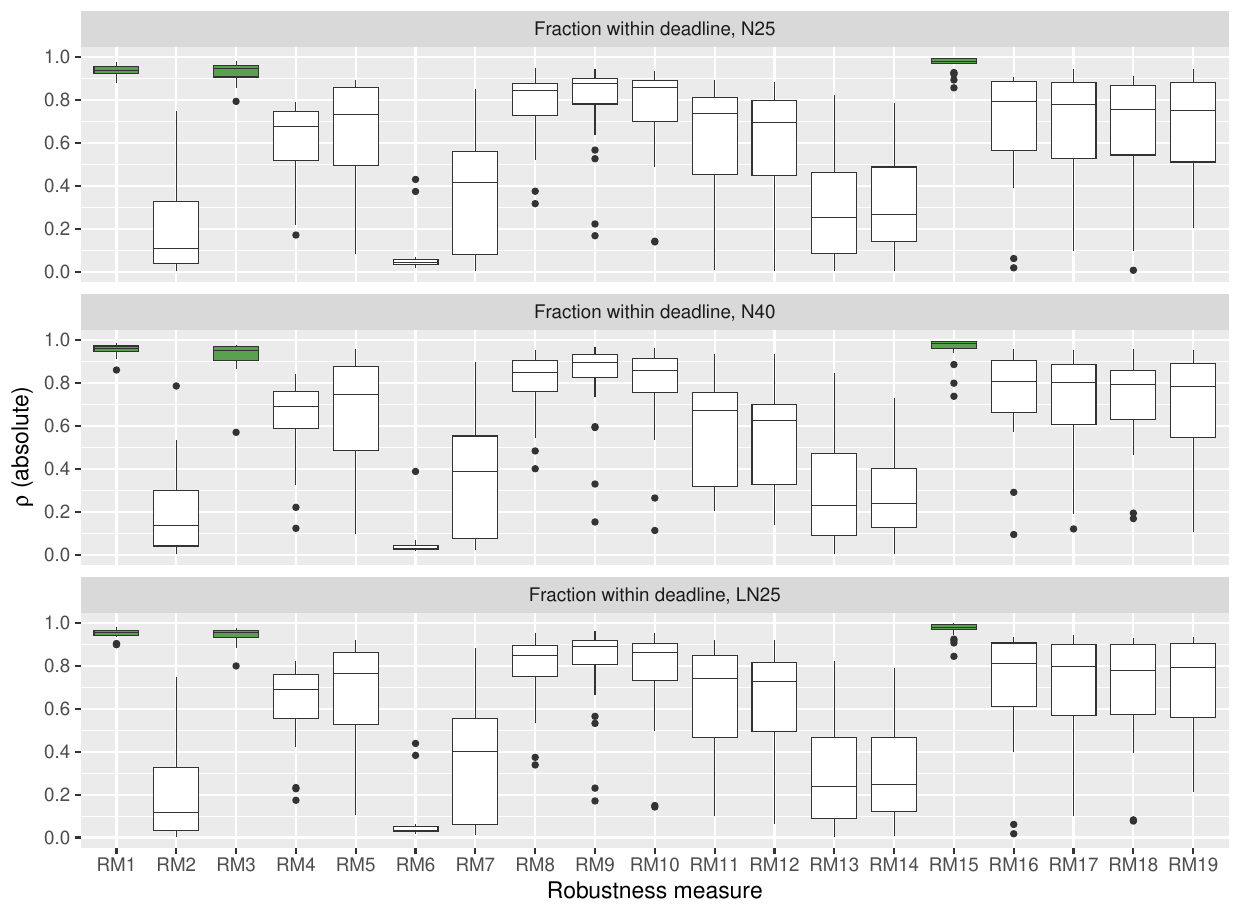}
    \caption{Spearman's coefficient $\rho$ (absolute value) for the average fraction within deadline. Green highlight indicates that the mean is greater than $0.9$.}
    \label{fig:corr_fraction_initial}
\end{figure}

Figures \ref{fig:corr_avgmakespan_initial} and \ref{fig:corr_fraction_initial} show the results for the average makespan and the fraction within deadline, respectively. We see that $RM_1$, $RM_3$, and $RM_{15}$ have the best overall performance on estimating both the average makespan and the fraction within deadline. $RM_6$, the sum of binary values indicating whether free slack exists, shows a correlation close to zero with both simulation measures and all distributions. This is because for most of the tested schedules, we inserted buffers for all jobs, giving each job at least some free slack. Only 10 out of 970 schedules per instance had no buffers inserted at all. This means that $RM_6$  only discerns between the 10 non-buffered schedules and the 960 others, resulting in a poor correlation.

The other measures show a high variability in correlation. Firstly, $RM_8$, $RM_9$ and $RM_{10}$ (weighted free slack with predecessors or successors) show a high median correlation with both simulation measures, but have some outlier instances where the correlation is very poor. Secondly,  $RM_2$ (sum of free slacks), $RM_7$ (weighted sum of free slacks with durations), and $RM_{13}$ and $RM_{14}$ (LP optimizing sum of buffer intervals) show low correlations. Lastly, $RM_4$, $RM_5$, $RM_{11}$, $RM_{12}$, $RM_{16}$, $RM_{17}$, $RM_{18}$ and $RM_{19}$ have a better median value but a low minimum value.

From these results, we see that the measures with the highest correlation with quality robustness are measures that either use the total slack or the approximated probability that the schedule finishes on time.

\subsection{Results solution robustness}\label{sec:first-exp-resul-sol}

\begin{figure}[htbp]
    \centering
    \includegraphics[width = \textwidth]{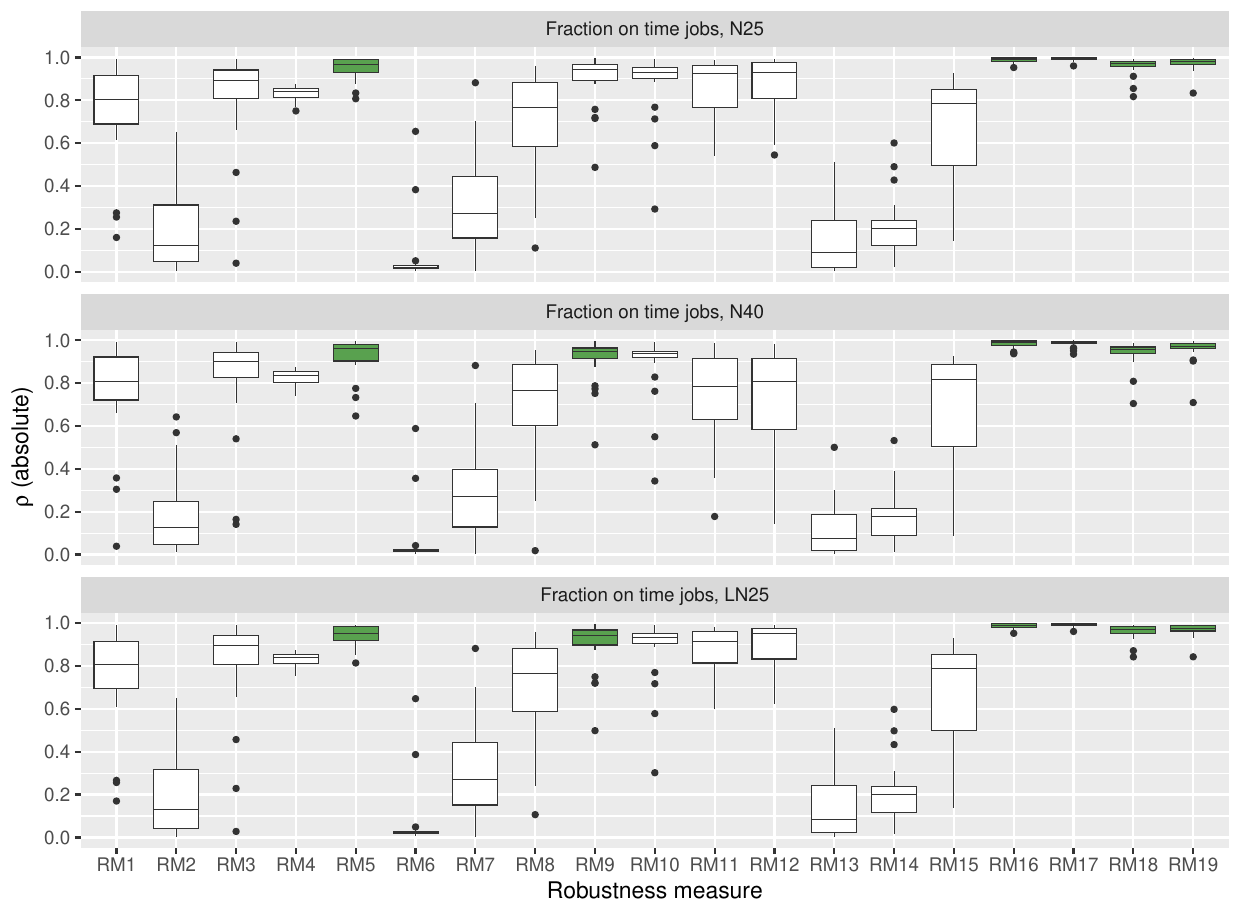}
    \caption{Spearman's coefficient $\rho$ (absolute value) for the average fraction of on-time jobs. Green highlight indicates that the mean is greater than $0.9$.}
    \label{fig:corr_onTime_initial}
\end{figure}

\begin{figure}[htbp]
    \centering
    \includegraphics[width = \textwidth]{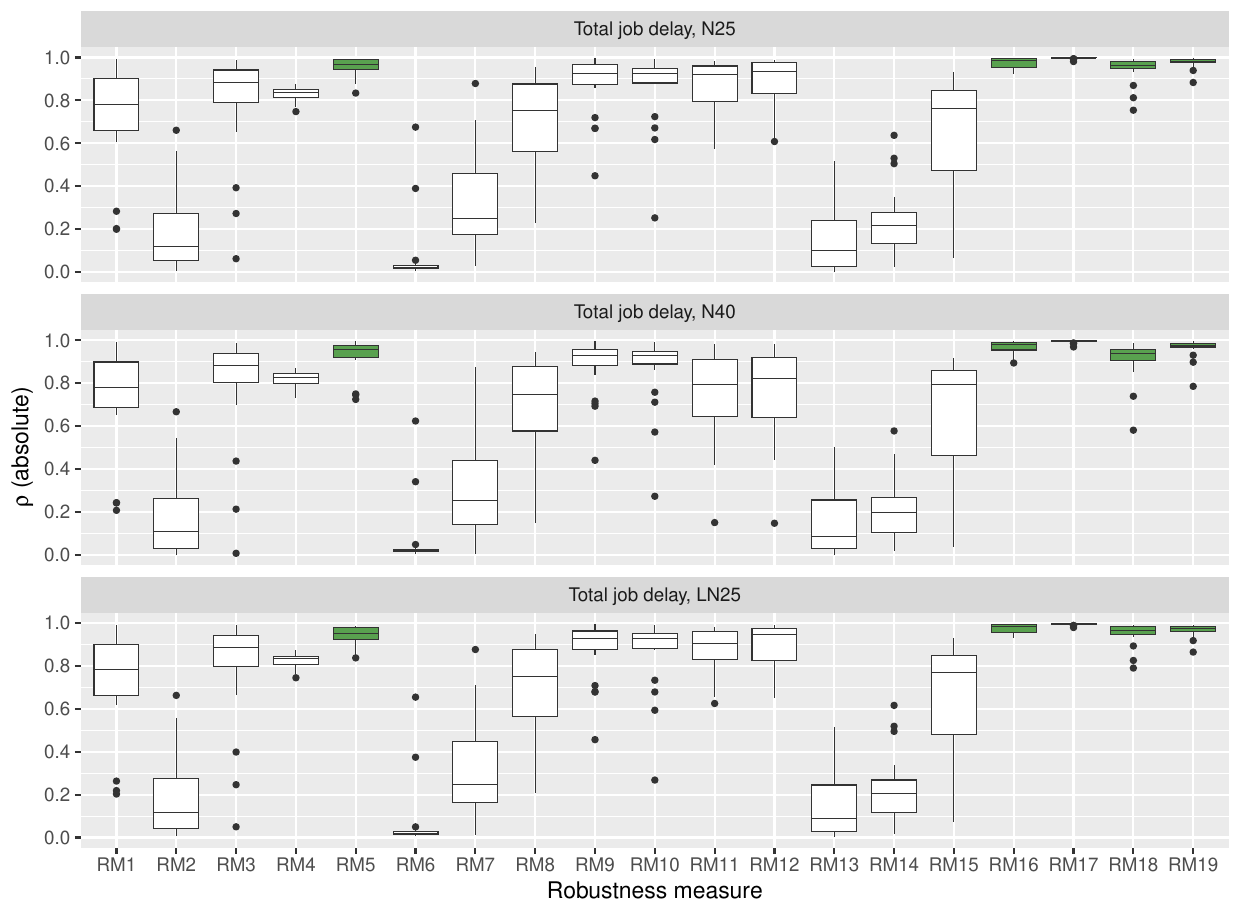}
    \caption{Spearman's coefficient $\rho$ (absolute value) for total job delay. Green highlight indicates that the mean is greater than $0.9$.}
    \label{fig:corr_totaldelay_initial}
\end{figure}

Figures \ref{fig:corr_onTime_initial} and \ref{fig:corr_totaldelay_initial} show the results for estimating the fraction of on-time jobs and the total job delay, respectively.

The figures indicate that $RM_5$, $RM_9$, $RM_{10}$, $RM_{16}$, $RM_{17}$, $RM_{18}$ and $RM_{19}$ have the best overall correlations with both the fraction of on-time jobs and the total job delay. Especially $RM_{16}$, $RM_{17}$, $RM_{18}$ and $RM_{19}$ show a very high correlation for both simulation measures. $RM_4$ shows a relatively high overall correlation, but the maximum value does not exceed $0.9$. $RM_{11}$ and $RM_{12}$ have high median correlation with both simulation measures when using $N_{25}$ and $LN_{25}$, but when using $N_{40}$ the minimal correlation is low. Again, $RM_6$ shows a correlation close to zero with both simulation measures and all distributions. The other measures show high variability in correlation for each combination of distribution and simulation measure. $RM_1$ and $RM_3$ show a high median correlation with both simulation measures, but have some outlier instances where the correlation is low. $RM_2$, $RM_7$, $RM_{13}$ and $RM_{14}$ show overall low correlations. $RM_8$ and $RM_{15}$ have a higher median value but a low minimum value.

From this, we see that multiple types of measures work well for solution robustness. Firstly, we have $RM_5$, $RM_{18}$, and $RM_{19}$ that consider the amount of free slack compared to the expected delay of a job. Then we have $RM_9$ and $RM_{10}$ that consider the free slack and number of successors of a job. Lastly, we have $RM_{16}$ and $RM_{17}$ that use normal approximation to approximate our simulation measures for solution robustness.

\subsection{Conclusion}

In this section, we investigated several robustness measures and their ability to predict the robustness of a schedule. We have shown that $RM_1$ (the sum of total slacks), $RM_3$ (minimum total slack), and $RM_{15}$ (normal approximation of on-time probability) all show high correlation with quality robustness. Of these three, the one based on normal approximation ($RM_{15}$) seems to perform the best. Furthermore, for solution robustness, we have shown that $RM_{16}$ (normal approximation of planned starting times), $RM_{17}$ (normal approximation for estimated starting delay), $RM_{18}$ (slack sufficiency of predecessors), and $RM_{19}$ (estimated starting delay) all show very high correlation. We also saw that $RM_5$ (sum of minimum free slack and expected delay), $RM_9$ (sum of free slack weighted with number of direct successors), and $RM_{10}$ (sum of free slack weighted with number of direct successors and processing time) show high correlation with solution robustness, though less so than the previous four. Lastly, we note that all well-performing measures are also competitive in terms of computation time. For all measures, evaluating 970 schedules of the biggest instance takes less than 100 milliseconds. We thereby also conclude that all these measures are applicable in a local search as they don't heavily impact the iteration speed. A full analysis of the computation times can be found in Appendix~\ref{sec:appendix-computation-time}.

\section{Local Search with Buffers}\label{sec:buffer-insertion}
Having identified the robustness measures that correlate well with both quality and solution robustness, we now turn to the question of how best to use these measures to create robust schedules. To this end, we present three local search methods using simulated annealing to create buffered schedules. One method uses an integrated approach to job assignment and buffer assignment, while the other two use a sequential approach. We will use the best measures from Section \ref{sec:Measure-correlation} as objectives for our local search and also compare them to find the best method for generating robust schedules. We will first describe our different local search algorithms.

\subsection{Local Search Neighbourhoods}

First, we describe the neighbourhoods used by our algorithms. Each job has an associated buffer value indicating its minimum required buffer time. With the machine assignments and the minimum buffer values, we can calculate all job starting times by passing through the schedule in topological order. So for our local search, we need to assign jobs to machines, determine the order in which they are executed, and assign buffer values to every job. To do this, we use 5 different local search neighbourhoods. We will use two neighbourhoods that modify the machine schedule assignment and machine ordering, and three that change the minimum job buffers. We end up with the following 5 neighbourhoods:

\begin{enumerate}
    \item $N_1$: Take two jobs and swap their positions.
    \item $N_2$: Take a random job, and insert it at a random position on a random machine.
    \item $N_3$: Decrease or increase the buffer of a random job.
    \item $N_4$: Take two jobs with a precedence relation between them, and move part of the buffer from one to the other.
    \item $N_5$: Take two jobs that follow each other directly on the same machine, and move part of the buffer from one to the other.
\end{enumerate}

The first two of these neighbourhoods are based on the $k$-move and $k$-swap moves often found in local search algorithms for Parallel Machine Scheduling (see, for instance, \citet{passage2025new}). For $N_3$, we choose a random job that is not the last on its machine and randomly increase or decrease its buffer by some random fraction of its processing time. We ensure that the buffer of a job is always capped at 0.7 times its processing time. For $N_4$ and $N_5$, we pick a random pair of jobs that either have a precedence relation or follow each other directly on the same machine, and then we randomly move part of one job's buffer to the other. This value is drawn from a normal distribution, with the mean set so that the two jobs end up with approximately equal buffers. As before, the buffer of a job remains capped at 0.7 times its processing time.

For the selection probabilities of our neighbourhoods, we use the ratio $4:4:3:2:2$ for $N_1:N_2:N_3:N_4:N_5$. This way, the relative frequency of neighbourhoods that alter job placement is roughly equal to that of neighbourhoods that alter the buffers. This is in accordance with \citet{lambrechts2008tabu}, who found that this results in the best solutions.

\subsection{Objective Function}

Next, we describe the objective function that we will use. Our objective function consists of two terms: an indicator function for quality robustness $qr$, and an indicator function for solution robustness $sr$. These two indicator values are calculated from our robustness measures, which we have normalized to fall between 0 and 1. These normalizations either use theoretical upper and lower bounds or are derived from experimental analysis of the ranges of the robustness measures. For example, $RM_{15}$ naturally falls between 0 and 1, so we do not have to normalize it. $RM_{16}$ is normalized by dividing by the number of jobs, as each summand falls between $0$ and $1$. And for $RM_{17}$, the best-case value is $0$, while the worst-case value can be anything. As we use $RM_{17}$ as a solution-robustness measure, we can then calculate the corresponding value of $sr$ as described in Equation~\eqref{eq:normalization} below:
\begin{equation}\label{eq:normalization}
    sr = \frac{\operatorname{initial}(RM_{17}) - RM_{17}(\sigma)}{\operatorname{initial}(RM_{17}) - 0}.
\end{equation}
Here, $\operatorname{initial}(RM_{17})$ refers to the initial value found for $RM_{17}$ when starting the local search procedure. As this initial value will be roughly equal to the worst value we find in the search procedure, this will set our normalized value $sr$ roughly between 0 and 1. After calculating the normalized values $qr$ and $sr$, we sum them using a relative weight $w$, giving the objective function in Equation~\eqref{eq:dynamic-objective-function}.

\begin{equation}\label{eq:dynamic-objective-function}
    f(qr, sr, w) = \frac{1}{w + 1}qr + \frac{w}{w + 1}sr.
\end{equation}

With this function, we can increase the importance of solution robustness by increasing $w$, nudging the local search algorithm towards more solution-robust schedules. We now have to determine the value of $w$ to use during our algorithm. One option is to use a fixed value for the entire local search, as done by \citet{liang2020integration}. However, initial experiments with different values of $w$ found that a fixed value is not well suited for our solution method. We have found that a high value of $w$ leads the algorithm to accept solutions with an excessively high makespan during the early, high-temperature phase of simulated annealing. Then, when the temperature decreases, it is unable to recover from these bad solutions. On the other hand, if we use a low value for $w$, we have found that the algorithm does not consider solution robustness enough, and hence the buffers cannot be optimized.

To address these problems, we use a dynamic value for $w$ that depends on the search progress, the current ratio of $qr$ to $sr$, and the number of precedence constraints. The dependence on search progress reflects our earlier observation: we want $w$ to start small and gradually increase as the search progresses. The dependence on the ratio $qr/sr$ reflects a second goal: when quality robustness is already high, we want the local search to focus more on solution robustness, and vice versa. Lastly, we found that for instances with many precedence relations compared to the number of jobs, quality robustness is more difficult to improve. Therefore, we want the objective to place more emphasis on $qr$ for these instances. Putting these considerations together, we arrive at the definition of $w$ given in Equation~\eqref{eq:obj-weight}. Here $i_{\text{current}}$ and $i_{\text{total}}$ indicate the current iteration and total number of iterations of the local search. Note that $w$ thus changes at every iteration. The variable $\alpha$ is an extra parameter that we can use to change the importance of solution robustness. In our experiments, we set $\alpha$ to 2.

\begin{equation}\label{eq:obj-weight}
    w = \qty(1 + \alpha\sqrt{\frac nr}\cdot\frac{i_{\text{current}}}{i_{\text{total}}})\cdot\frac{qr}{sr}
\end{equation}

\subsection{Local Search Algorithms}

Now that we have explained the neighbourhoods and objective function we will use, we can continue with our different local search algorithms. We start with an algorithm that simultaneously uses job placement and buffer neighbourhoods (CSA). For this algorithm, we start by applying a greedy algorithm to find an initial solution. We order the jobs by earliest possible start time based on release dates and precedence constraints, and then assign them in that order to the machine with the lowest makespan. The buffers of the jobs are initialised at 0.1 times their average processing time. This greedy algorithm usually produces a schedule with a deterministic makespan that is below the deadline. If this is not the case, we use a hill-climbing procedure that moves single jobs until the makespan is below the deadline. After this, we start with the simulated annealing algorithm. First, we calculate the solution weight corresponding to the current solution and iteration. Then we apply a randomly selected neighbour to our solution. We calculate the new robustness scores for this altered schedule and calculate the difference with the previous solution. If the new solution is better, we always accept it. If not, we accept it with a probability that is based on the current temperature and the size of the difference. Every time the current iteration is divisible by our cooling interval, we multiply the current temperature by our cooling rate. We stop when the temperature is lower than the given minimum. At that point, we return the best solution that we have encountered during the algorithm. A full overview of this algorithm can also be found in Algorithm~\ref{alg:CSA}.

\begin{algorithm}
\caption{Local Search with Combined Simulated Annealing}
\label{alg:CSA}
\begin{algorithmic}[1]
\Require Problem instance $I$, weight parameter $\alpha$, initial temperature $T_0$, minimum temperature $T_{\min}$, cooling rate $CR$, cooling interval $CI$

\Procedure{LocalSearch}{$I, \alpha, T_0, T_{\min}, CR, CI$}
    \State $S \leftarrow \Call{GreedyInitialSolution}{I}$
    \If{$\operatorname{Makespan}(S) > d$}
        \State $S \leftarrow \Call{MakespanHillclimber}{S}$
    \EndIf
    \State $(qr, sr) \leftarrow \Call{CalculateRobustness}{S}$ \Comment{Start of simulated annealing}
    \State $i \leftarrow 0$
    \State $T \leftarrow T_0$
    \State $i_\text{total} = \log_{CR}(T_{\text{min}} / T_0) * CI$
    \While{$T > T_{\min}$}
        \State $i \leftarrow i + 1$
        \State $w \leftarrow \Call{GetSolutionWeight}{qr, sr, i, i_{\text{total}}, \alpha}$
        \State $S' \leftarrow \Call{ChooseAndApplyNeighbor}{S}$
        \State $(qr', sr') \leftarrow \Call{CalculateRobustness}{S'}$
        \State $\Delta \leftarrow f(qr', sr', w) - f(qr, sr, w)$
        \State $p_{accept} \leftarrow \Call{Random}{0, 1}$
        \If{$\Delta \geq 0 \lor e^{-\left(\frac{\Delta}{T}\right)} < p_{accept}$}
            \State $(S, qr, sr) \leftarrow (S', qr', sr')$
        \EndIf
        \If{$i \mod CI = 0$}
            \State $T \leftarrow CR \cdot T$
        \EndIf
    \EndWhile
    \State \Return best encountered solution $S^*$
\EndProcedure
\end{algorithmic}
\end{algorithm}

Our other two local search algorithms have separate phases for solving the job placement and assigning job buffers. The first of these two is referred to as Separated Simulated Annealing (2SA). From our initial solution, we start with a round of simulated annealing using only $N_1$ and $N_2$, both occurring with equal probability. In this round of simulated annealing, we also set the solution robustness weight to 0, thus only optimizing for quality robustness. Then we have a second round of simulated annealing. In this second round, we switch to the buffer neighbourhoods $N_3$, $N_4$, and $N_5$, and we use the same objective as in the combined simulated annealing. Both of these rounds of simulated annealing use the same cooling schedule as the combined simulated annealing algorithm, but since it has two simulated annealing rounds, the algorithm will take twice as long. The idea behind this method is that by starting with just the job placement and focusing only on quality robustness, we get a good baseline schedule into which we can then insert the buffers. An overview of this algorithm can also be found in Algorithm~\ref{alg:2SA}.

\begin{algorithm}
\caption{Local Search with Separated Simulated Annealing}
\label{alg:2SA}
\begin{algorithmic}[1]
\Require Problem instance $I$, weight parameter $\alpha$, initial temperature $T_0$, minimum
temperature $T_{\min}$, cooling rate $CR$, cooling interval $CI$

\Procedure{LocalSearchSeparated}{$I, \alpha, T_0, T_{\min}, CR, CI$}
    \State $S \leftarrow \Call{GreedyInitialSolution}{I}$
    \If{$\operatorname{Makespan}(S) > d$}
        \State $S \leftarrow \Call{MakespanHillclimber}{S}$
    \EndIf    \State $\Call{SetNeighbourhoodOperators}{\{N_1: 1, N_2: 1\}}$
    \State $S \leftarrow \Call{SimulatedAnnealing}{S, T_0, T_{\min}, CR, CI, w = 0}$
    \State $\Call{SetNeighbourhoodOperators}{\{N_3: 3, N_4: 2, N_5: 2\}}$    
    \State $S \leftarrow \Call{SimulatedAnnealing}{S, T_0, T_{\min}, CR, CI}$
    \State \Return best encountered solution $S^*$
\EndProcedure
\end{algorithmic}
\end{algorithm}

The last algorithm we consider is an iterated local search algorithm (ILS). In this algorithm, we still have separate phases for job placement and buffer distribution, but we keep switching between them. After we have initialised the schedule, we start with a simulated annealing run only using $N_1$ and $N_2$. Then we do a simulated annealing run for just the buffer insertion. After this, we perturb the solution by applying several random neighbourhoods from the whole set $N_1 - N_5$ to help the solution escape local optima. This process is repeated 8 times. For the first run, we use starting temperature $T_0$, but for the repeats, we use a lower starting temperature $T_r$, which we set as the geometric mean of the original starting temperature and the minimum temperature, $\sqrt{T_0\cdot T_{\min}}$. We also do not fix $w$ to 0 during the job placement phase, since we repeatedly switch between the two neighbourhood types, meaning solution robustness must be considered throughout. The idea is that by iterating the two operation types, we can account for interactions between changing job placement and changing the buffers. A short overview of this algorithm can be found in Algorithm~\ref{alg:ILS}.

\begin{algorithm}
\caption{Iterated Local Search with Separated Simulated Annealing}
\label{alg:ILS}
\begin{algorithmic}[1]
\Require Problem instance $I$, weight parameter $\alpha$, initial temperature $T_0$, minimum
temperature $T_{\min}$, restart temperature $T_r$ cooling rate $CR$, cooling interval $CI$,
\Procedure{IteratedLocalSearch}{$I, \alpha, T_0, T_{\min}, CR, CI$}
    \State $S \leftarrow \Call{GreedyInitialSolution}{I}$
    \If{$\operatorname{Makespan}(S) > d$}
        \State $S \leftarrow \Call{MakespanHillclimber}{S}$
    \EndIf    \State $n \leftarrow 0$
    \State $T \leftarrow T_0$
    \While{$n < 8$}
        \State $n \leftarrow n + 1$
        \State $\Call{SetNeighbourhoodOperators}{\{N_1: 1, N_2: 1\}}$
        \State $S \leftarrow \Call{SimulatedAnnealing}{S, T, T_{\min}, CR, CI}$
        \State $\Call{SetNeighbourhoodOperators}{\{N_3: 3, N_4: 2, N_5: 2\}}$
        \State $S \leftarrow \Call{SimulatedAnnealing}{S, T, T_{\min}, CR, CI}$
        \State $S \leftarrow \Call{PerturbSolution}{S}$
        \State $T \leftarrow T_r$
    \EndWhile
    \State \Return best encountered solution $S^*$
\EndProcedure
\end{algorithmic}
\end{algorithm}

\subsection{Experimental Results}
To decide on the best way to create schedules that are both quality robust and solution robust, we will use our three local search algorithms to generate schedules for all instances from Section~\ref{sec:Measure-correlation}. We use $RM_1$, $RM_3$ and $RM_{15}$ for the quality robustness measures, and we use $RM_5$, $RM_9$, $RM_{10}$, $RM_{16}$, $RM_{17}$, $RM_{18}$ and $RM_{19}$ for our solution robustness measures. With these, we will investigate which combination of algorithm type, quality measure, and solution measure generates the best schedules. For each algorithm, we set the initial temperature to 0.01 and the minimum temperature to $10^{-5}$. Our cooling rate will be 0.98. For the Combined Simulated Annealing and the Separated Simulated Annealing, we use a cooling interval of 200 iterations; both algorithms run each simulated annealing phase for a full cooling schedule from $T_0$ to $T_{\min}$. For the Iterated Local Search, however, we shorten the interval to 50 iterations, since here the cooling schedule is repeated 8 times within the same overall iteration budget, and a 200-iteration interval per repetition would give ILS a much larger total number of iterations than the other two algorithms. Although the total number of iterations for each algorithm still differs, they are comparable, ranging between 75,000 and 150,000.

\begin{figure}
    \centering
    \includegraphics[width = \linewidth]{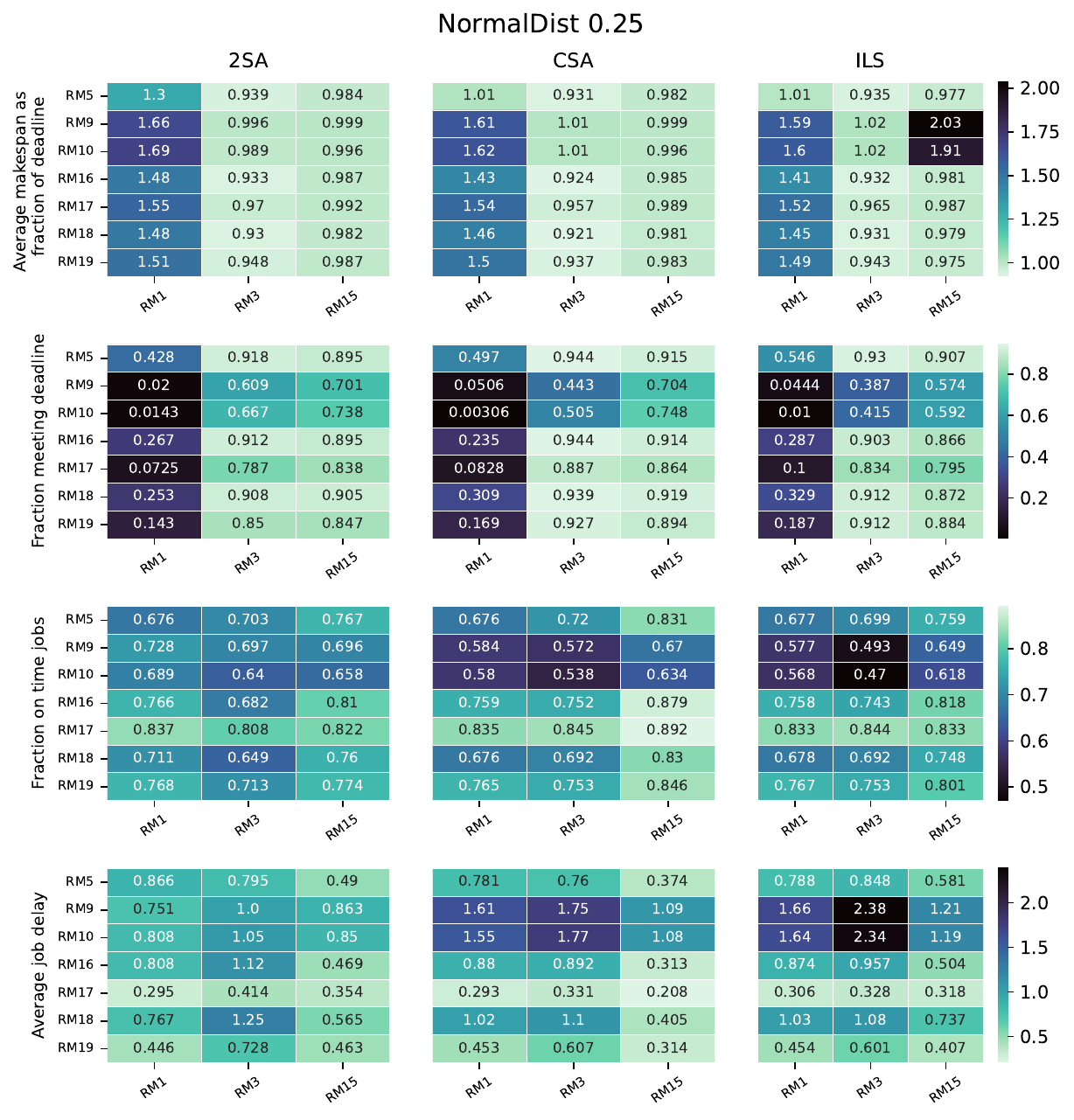}
    \caption{Results from the local search experiment, where every number corresponds to a simulation average for a specific combination of algorithm type and objective function. The distributions used for the processing times in this figure are all normal distributions with a standard deviation of 0.25 times the average duration. A lighter colour indicates a more robust score, while darker indicates a less robust score.}
    \label{fig:results-ls-norm25}
\end{figure}

\begin{figure}
    \centering
    \includegraphics[width = \linewidth]{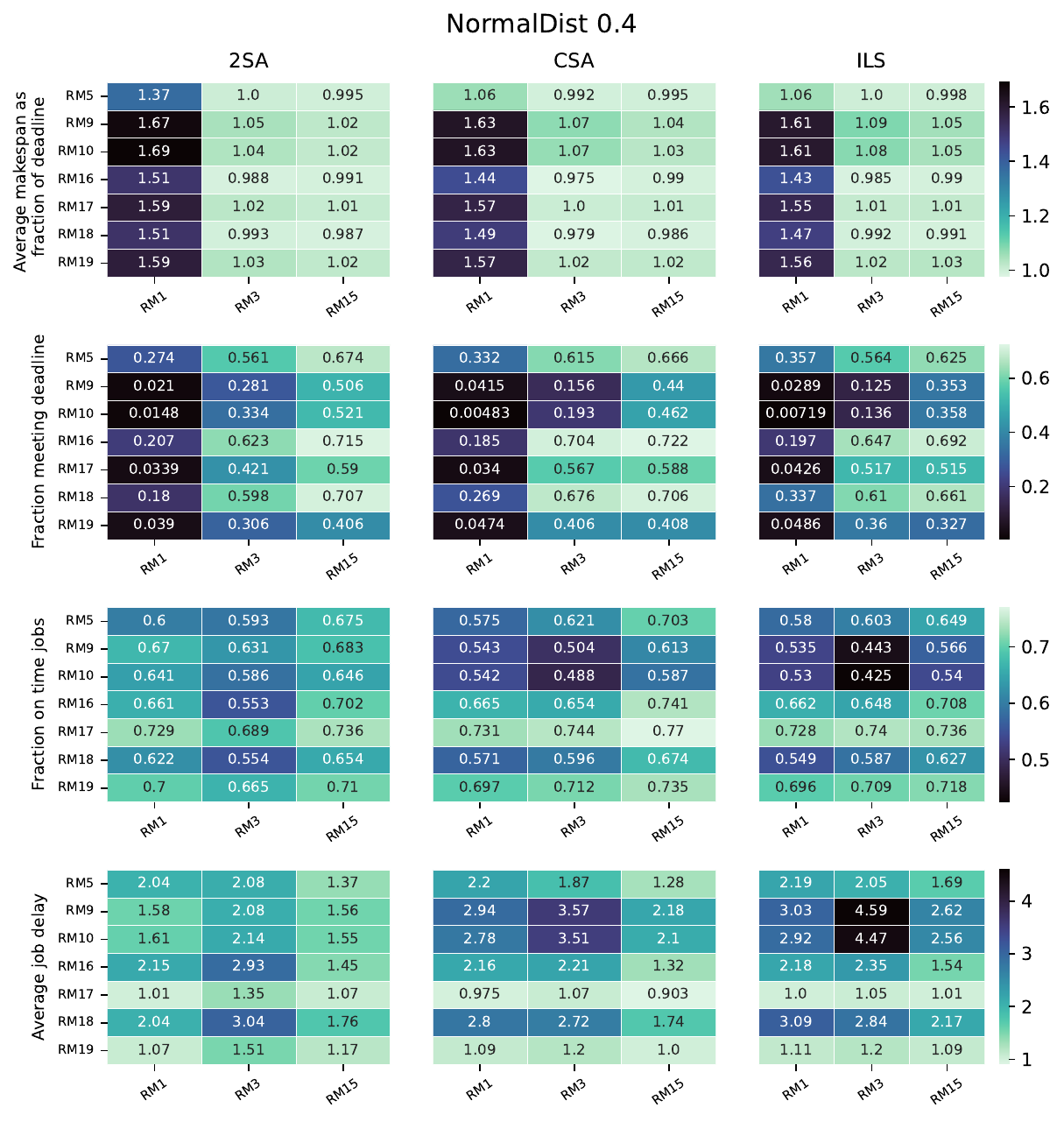}
    \caption{Results from the local search experiment, where every number corresponds to a simulation average for a specific combination of algorithm type and objective function. The distributions used for the processing times in this figure are all normal distributions with a standard deviation of 0.4 times the average duration. A lighter colour indicates a more robust score, while darker indicates a less robust score.}
    \label{fig:results-ls-norm40}
\end{figure}

\begin{figure}
    \centering
    \includegraphics[width = \linewidth]{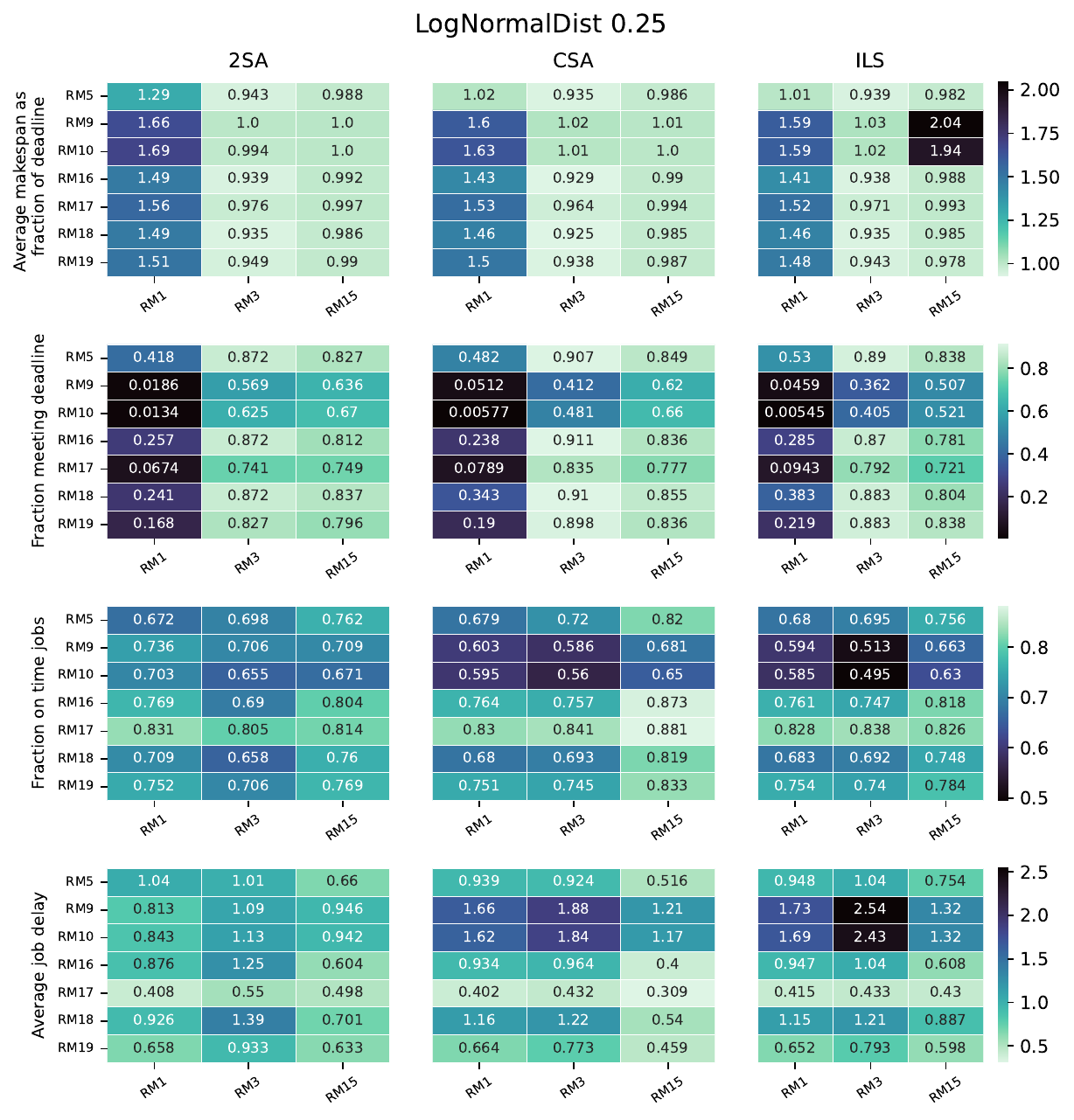}
    \caption{Results from local search experiment, where every number corresponds to a simulation average for a specific combination of algorithm type and objective function. The distributions used for the processing times in this figure are all log-normal distributions with a standard deviation of 0.25 times the average duration. A lighter colour indicates a more robust score, while darker indicates a less robust score.}
    \label{fig:results-ls-lognorm25}
\end{figure}

For each instance, we run each local search algorithm with every possible combination of robustness measures for the objectives giving a total of 21 combinations of measures. We repeat this 20 times for every instance. We then simulate all generated schedules 10,000 times and note the average makespan, average fraction within deadline, average fraction of on-time jobs, and the average total delay. These simulations use the same execution policy as in Section~\ref{sec:Measure-correlation}, where jobs are not allowed to begin before their planned starting time. We report the results of these simulations in Figure~\ref{fig:results-ls-norm25}, Figure~\ref{fig:results-ls-norm40} and Figure~\ref{fig:results-ls-lognorm25}. Each figure corresponds to a different distribution used during simulation. Every figure consists of a set of twelve larger blocks. Each block represents the results for one simulation measure (spread over the rows) and one of the local search algorithms (spread over the columns). In each large block, there are 21 smaller blocks corresponding to a specific combination of robustness measures in the objectives. For these measures, the rows correspond to the different solution measures, while the columns correspond to the different quality measures. Each number corresponds to the average taken over all runs of the corresponding algorithm for all instances. To normalize the data, we have divided the average makespan by the deadline of the corresponding instance.

Let us first compare the different algorithms. We first compare the Iterated Local Search algorithm with the Combined Simulated Annealing algorithm. Looking at the averages of all individual measure combinations, CSA almost always scores better. Note that for the fraction meeting deadline and fraction on time jobs, a higher score is better, while for the other two, a lower score is better. When ILS does score better on average, it only does so marginally. Thus, we can conclude that the CSA algorithm is better than the ILS algorithm. We are left with comparing the Combined Simulated Annealing algorithm and the Separated Simulated Annealing algorithm. For these two, the difference between combinations of measures appears somewhat larger, as CSA performs clearly better for some measures, while 2SA is better for others. However, if we disregard the worst-performing measures overall ($RM_1$ for quality and $RM_9$ and $RM_{10}$ for solution), CSA outperforms separated simulated annealing in both quality and solution robustness across all remaining measure combinations. 

Next, we compare the different combinations of robustness measures. We first note that $RM_1$ performs very badly for the quality robustness measures. For both the 2SA and CSA algorithms, the quality robustness of the generated schedules is very low, while solution robustness is somewhat comparable to the other measures. This makes $RM_1$ clearly inferior to $RM_3$ and $RM_{15}$. We can also immediately see that $RM_9$ and $RM_{10}$ perform worse than the other solution measures in both quality and solution robustness. This leaves us with $RM_3$ and $RM_{15}$ for quality robustness, and $RM_5$, $RM_{16}$, $RM_{17}$, $RM_{18}$, and $RM_{19}$ for solution robustness.

Between $RM_3$ and $RM_{15}$, there is a slight difference between the different processing time distributions. Although $RM_{15}$ mostly outperforms $RM_3$ for the $N_{40}$ distribution, $RM_3$ produces more quality robust schedules for $N_{25}$ and $LN_{25}$. It is therefore difficult to determine which of these two is the best overall robustness measure. For the solution robustness measures, the largest differences appear for the distribution $N_{40}$ in Figure~\ref{fig:results-ls-norm40}. Here, the schedules have the highest quality robustness with $RM_{16}$ as the solution robustness measure, while the solution robustness is higher with $RM_{17}$. The remaining solution robustness measures ($RM_5$, $RM_{18}$ and $RM_{19}$) are all completely inferior to either one of these. Between $RM_{16}$ and $RM_{17}$, it is again difficult to say which of these is the absolute best, as we are in a similar situation as between $RM_3$ and $RM_{15}$

\begin{figure}
    \centering
    \includegraphics[width=\linewidth]{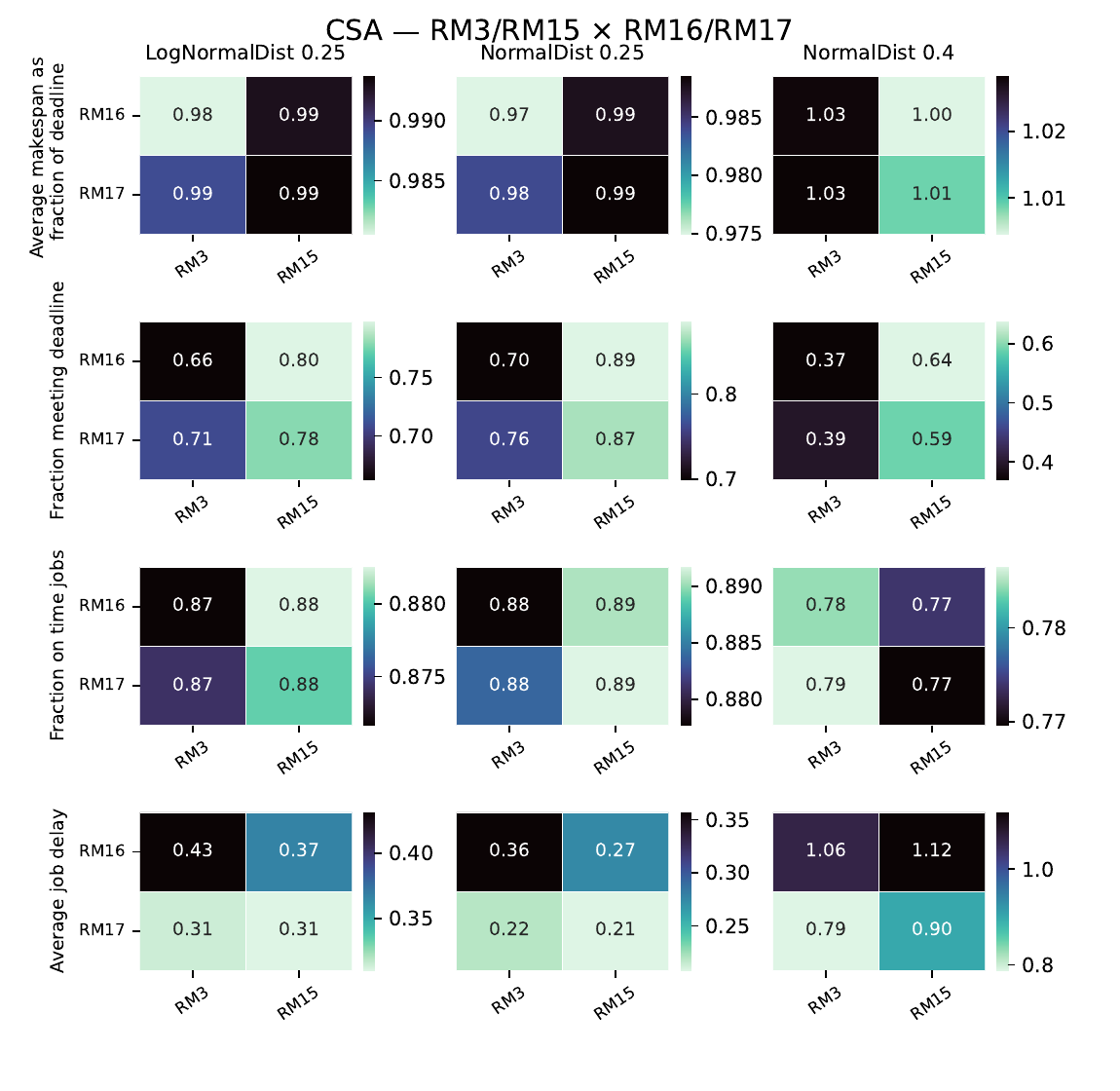}
    \caption{Results from second local search experiment where we compared $RM_3$ and $RM_{15}$ as quality robustness measures, and $RM_{16}$ and $RM_{17}$ as solution robustness measures. Every number corresponds to a simulation average for a specific combination of processing time distribution type and objective function. A lighter colour indicates a more robust score, while darker indicates a less robust score.}
    \label{fig:results-best-measures}
\end{figure}

To discern between these final sets of robustness measures, we perform a final experiment where we perform a similar comparison as before. However, we make a small change in how we set it up. We increase the $\alpha$ parameter in the objective weight function when $RM_{3}$ or $RM_{16}$ is used. Because these had worse results for solution robustness but better results for quality robustness, we adjust $\alpha$ only for these measures to better compare $RM_{15}$ and $RM_{17}$. For this experiment, we only use the CSA algorithm, as it was found to be the best algorithm in the previous experiment. Figure~\ref{fig:results-best-measures} shows the results of this experiment.

If we first compare $RM_{3}$ and $RM_{15}$, we see a clear difference. Schedules produced using $RM_{15}$ are always clearly more quality-robust and either also more solution-robust, or at worst, roughly equally solution-robust as the schedules produced using $RM_3$. Even when the results for $RM_{15}$ are better, the difference is negligibly small. Thus we find that $RM_{15}$ is the best quality robustness measure. Notably, this holds despite the increased $\alpha$ value used for $RM_3$, which was intended to favour its solution robustness performance, thus reinforcing that $RM_{15}$'s advantage is not an artifact of this adjustment. For $RM_{16}$ and $RM_{17}$, it is more difficult to tell the difference. For most simulation measures, they are almost equal. The only measure with a clear difference is the average job delay, with $RM_{17}$ performing slightly better. Overall, however, $RM_{16}$ and $RM_{17}$ seem to perform almost identically. This is not unexpected, as both measures are based on the same type of approximation. The only difference is that one of them uses the expected value, while the other uses the probability. These tests do not show a noticeable difference in performance. However, we can conclude that normal approximation gives us the best results for both quality and solution robustness when using normal and lognormal distributions for the job processing times.

\section{Conclusion}\label{sec:conc}
In this paper, we addressed the problem of creating robust schedules for the Stochastic Parallel Machine Scheduling Problem (SPMSP). Given the many robustness measures and scheduling strategies proposed in the literature, with no clear consensus on which perform best, we set out to do two things. Firstly, we compared new and existing robustness measures on both quality and solution robustness to identify which measures work well within a local search algorithm. Secondly, we compared different strategies for distributing slack within a local search algorithm for the SPMSP.

We started with an extensive simulation study where we evaluated 13 existing and 6 new robustness measures, comparing their correlation with both quality robustness and solution robustness. From this evaluation, we found 3 measures highly correlated with quality robustness and 7 that were highly correlated with solution robustness. All of these measures were also easy to compute, requiring a limited number of traversals through the schedule for each measure. To the best of our knowledge, this is the first time such a large comparison has been done for both quality robustness and solution robustness.

Next, we took all the highly correlated measures and tested them as objectives in various local search algorithms to find the best way to create schedules that are both quality-robust and solution-robust. We introduced a local search procedure with dynamic solution weights to balance the quality and solution robustness. We implemented this in three local search algorithms to compare integrated and sequential methods for assigning jobs and buffers. We found that the Combined Simulated Annealing (CSA) algorithm, which simultaneously assigns buffers and jobs, created the most robust schedules. The dynamic solution weight we introduced allowed the search to adaptively balance both robustness types over the course of the search. Furthermore, we found that most robust schedules were generated when using $RM_{15}$ (approximation of probability that the schedule finishes on time) as a quality robustness measure and using either $RM_{16}$ (sum of approximated probabilities that jobs start on time) or $RM_{17}$ (sum of approximated expected starting delays of jobs) as solution measures. All three measures are based on an approximation for the maximum of two normal distributions as in \citet{passage2025new}, indicating that this is a very effective method for measuring the robustness of a schedule.

For future work, it would be interesting to explore why some measures are more effective than others in generating robust schedules. This would make it easier to predict how a general measure function might affect a scheduling algorithm. Moreover, it would also be interesting to see how these results compare to different problems. In particular, as the SPMSP is an abstract version of many practical scheduling problems, it would be interesting to see whether these results carry over to real-world applications. Furthermore, it would also be interesting to see if the results change for more complex theoretical scheduling problems.

\section*{Acknowledgements}
This work is part of the NWO LTP-ROBUST RAIL Lab, a collaboration between Utrecht University, the Delft University of Technology, NS, and ProRail. More information at \url{https://icai.ai/lab/rail-lab-utrecht/}.

\appendix

\section{Details on Instance and Schedule Generation}
\subsection{Generating the instances} \label{sec:app-instance}
For these experiments, we have used a total of 24 instances of the SPMSP. We first took the 12 instances from \citet{passage2025new}, consisting of 6 small instances with $n = 30$ jobs and 6 large instances with $n = 100$ jobs. The small instances have further parameters $r\in \{15, 30, 75\}$ for the number of precedence constraints and $m\in\{4, 8\}$ for the number of machines. For the large instances we have $r\in\{50, 100, 250\}$ and $m\in\{6, 12\}$. For each combination of parameters $n,m,r$, we also generated an additional instance. Average processing times for all jobs, $p_j$, are chosen uniformly at random from the interval $[1, 20]$, and release dates are chosen uniformly at random from the interval $[0, \floor{n/2}]$. At this point, we do not yet specify the processing time distributions. We start by generating schedules using the averages $p_j$, and then consider several distributions for simulation. Precedence constraints are assigned randomly, ensuring that no cycles form. Deadlines are generated by using the two methods of \citet{van2006trade}, and then taking the maximum of the two obtained values. The first method uses the fact that the minimum length of the schedule equals the sum of the processing times divided by the number of machines. But we also account for the release dates by adding the $m$ smallest release dates to the sum of the processing times:
\begin{equation}
    \ell_{min} = \frac{\sum_{j = 1}^np_j + \sum_{j = 1}^mr_j}{m}.
\end{equation}
Here we assume that the jobs are numbered such that $r_1\leq r_2\leq\dots\leq r_n$. \citet{van2006trade} found that, to obtain a stable schedule, we need to multiply this by 1.3 to account for uncertainty in processing times. The second way to get a lower bound for the deadline is to find the length of the critical path $\ell_{cp}$. This is the longest path of activities in the instance, where every activity must be completed before the next can start. To account for the stochastic processing times, \citet{van2006trade} found that it is necessary to scale the length of this path by a factor of 1.5. However, since the standard deviation of the critical path also depends on the number of jobs on the path, we choose to only increase this length by $(0.5\cdot \ell_{cp})/ \sqrt{n_{cp}}$, where $n_{cp}$ is the number of jobs in the critical path. Then we set the deadline for the instance as:
\begin{equation}
    d = \max\Big\{\ell_{cp}\cdot\qty(1 + \frac{0.5}{\sqrt{n_{cp}}}), \ell_{min}\cdot 1.3\Big\}.
\end{equation}

\subsection{Generating the schedules} \label{sec:app-scehdules}

To accurately test our robustness measures, we need a wide range of different schedules for each instance. An initial solution is obtained using a greedy algorithm that repeatedly selects the job that can start earliest from all jobs that have no unassigned predecessors. This selected job is then scheduled on the machine with the earliest completion time, ties are broken randomly. Next, we use a hill-climbing procedure which considers neighbour solutions obtained by either swapping a job with another job on the same machine ($N_1$) or by moving a job to any feasible position on any machine ($N_2$). When exploring a neighbourhood, the neighbours are selected in a random order, and the first improving neighbour is then accepted for the new schedule. The two neighbourhood types are explored in an alternating fashion, starting with $N_1$. The resulting schedule is returned when no more improvements can be found. If the deterministic makespan is larger than the deadline, the procedure is restarted.

This local search procedure always produces earliest start schedules, as the algorithm aims to minimize the makespan. However, we want to have schedules with buffers for every job. Therefore, we have to alter the starting times of all jobs to introduce buffers into the schedule. To achieve this, we make use of the LP from $RM_{14}$. This LP assigns each job $j$ ($j = 1,\dots, n$) an interval $[e_j, l_j]$ of possible start times. If we then set the starting time of job $j$ to $e_j$, we have given each job a buffer of size $l_j - e_j$. The objective ensures that these buffers are evenly distributed between the different jobs. By solving this LP, we always obtain the same buffers given the same earliest start schedule. 

To obtain a variety of earliest start schedules, we proceed as follows: for every job $j$ ($j = 1,\dots, n$) in the schedule $\overline\sigma$ we just obtained, we now have a starting time $ s_j^{\overline\sigma}$, a processing time $p_j$, and a buffer time $b_j^{\overline\sigma}$. Then, to create a schedule with more diverse buffers, we multiply each job's buffer by a factor $\mu_j < 1$ drawn uniformly at random from a predetermined interval $[a, b]$. Now each job has an associated minimum buffer time between the end of that job and the starting times of its successor jobs. Using these minimum buffer times, we pass through the schedule in topological order and, for each job, compute the earliest starting time that respects its predecessors' buffers. We repeat this 5 times with a total of 19 intervals $[a,b]$. We start with the interval $[0, 0.1]$. Then we incrementally increase the upper bound by $0.1$, giving us the intervals $[0, 0.2], [0,0.3], \dots [0,1]$. After this, we start increasing the lower bound by $0.1$ while keeping the upper bound at $1$ giving us the intervals $[0.1, 1], [0.2, 1],\dots, [0.9,1]$. We also include the schedule with maximum buffers and the schedule with 0 buffers. This gives us a total of $19\times 5 + 2 = 97$ buffer schedules per earliest start schedule. Although this method relies on one of the robustness measures ($RM_{14}$) for generating the buffer intervals, the subsequent per-job independent randomization substantially perturbs the relative buffer distribution across jobs. As this is precisely the objective that $RM_{14}$ optimizes for, we believe this mitigates the risk of artificially inflating $RM_{14}$'s performance. Notably, our results in Sections~\ref{sec:first-exp-results-quality} and \ref{sec:first-exp-resul-sol} show that $RM_{13}$ and $RM_{14}$ in fact perform comparatively poorly, which is inconsistent with a strong circularity effect.

\section{Computation Time}\label{sec:appendix-computation-time}
To use a robustness measure as objective in a local search algorithm, the measure must be quickly computable. Longer computation time means more time per iteration and hence a slower local search algorithm. To compare the computational efficiency of the robustness measures, we record the total time needed to evaluate the 970 generated schedules for the instances $30j-15r-4m$ and $100j-250r-12m$. We repeat this 100 times and report the average total time. For comparison, we also record the time needed to evaluate the 970 schedules by performing 100 simulation runs. For this experiment, a laptop with a 12th Gen Intel\textsuperscript{\textregistered} Core\textsuperscript{TM} i7-1265U @ 1.80GHz processor was used.

Results are shown in Table \ref{tab:times}. They show that the basic slack-based measures $RM_1$ to $RM_{10}$ have similar computational times. To compute them, the schedule has to be traversed once. $RM_{18}$ and $RM_{19}$ are slack-based measures that take more time because the schedule is traversed twice. Additionally, we can see that normal approximation methods $RM_{15}$, $RM_{16}$, and $RM_{17}$ are less computationally efficient than the simpler slack-based measures. This is likely due to the computational overhead from calculating the distributions. However, they differ only by a factor of 3 to 4, making them still feasible to use in a local search algorithm. The more complex slack-based measures $RM_{11}$ and $RM_{12}$, as well as the LP-based measures $RM_{13}$ and $RM_{14}$, have poor efficiency compared to the other measures. Finally, all robustness measures except for $RM_{13}$ and $RM_{14}$ show significantly better computational efficiency than performing 100 simulations. This is as expected, since running 100 simulations requires traversal through the schedule 100 times. For $RM_{13}$ and $RM_{14}$, it is not unexpected that the performance is so much worse, as a complete LP needs to be solved each iteration. It may be possible to improve the performance of this measure by reusing old LP solutions, but we have not tested this.

\begin{table}[htbp]
    \centering
    \caption{Computation time in milliseconds of evaluating 970 schedules}
    \begin{tabular}{@{}lrr}
        \toprule
         & \textsc{30j-15r-4m} & \textsc{100j-250r-12m} \\ \midrule
        \multicolumn{1}{l}{$RM_1$}       & 4.17    & 29.99    \\
        \multicolumn{1}{l}{$RM_2$}       & 3.21    & 24.55    \\
        \multicolumn{1}{l}{$RM_3$}       & 4.20    & 30.07    \\
        \multicolumn{1}{l}{$RM_4$}       & 3.12    & 24.57    \\
        \multicolumn{1}{l}{$RM_5$}       & 4.38    & 28.91    \\
        \multicolumn{1}{l}{$RM_6$}       & 3.09    & 24.40    \\
        \multicolumn{1}{l}{$RM_7$}       & 3.09    & 24.49    \\
        \multicolumn{1}{l}{$RM_8$}       & 3.12    & 24.60    \\
        \multicolumn{1}{l}{$RM_9$}       & 3.11    & 24.50    \\
        \multicolumn{1}{l}{$RM_{10}$}    & 3.08    & 24.69    \\
        \multicolumn{1}{l}{$RM_{11}$}    & 19.81   & 267.61   \\
        \multicolumn{1}{l}{$RM_{12}$}    & 19.66   & 268.65   \\
        \multicolumn{1}{l}{$RM_{13}$}    & 401.70  & 2171.15  \\
        \multicolumn{1}{l}{$RM_{14}$}    & 650.84  & 2611.35  \\
        \multicolumn{1}{l}{$RM_{15}$}    & 11.78   & 84.85    \\
        \multicolumn{1}{l}{$RM_{16}$}    & 13.68   & 96.01    \\
        \multicolumn{1}{l}{$RM_{17}$}    & 12.98   & 92.48    \\
        \multicolumn{1}{l}{$RM_{18}$}    & 5.31    & 43.11    \\
        \multicolumn{1}{l}{$RM_{19}$}    & 6.05    & 46.31    \\
        \multicolumn{1}{l}{100Sim}       & 327.33  & 1463.30  \\ \bottomrule
    \end{tabular}
    \label{tab:times}
\end{table}

\bibliographystyle{elsarticle-harv}
\bibliography{refs}

@article{chtourou2008two,
  title={A two-stage-priority-rule-based algorithm for robust resource-constrained project scheduling},
  author={Chtourou, H{\'e}di and Haouari, Mohamed},
  journal={Computers \& industrial engineering},
  volume={55},
  number={1},
  pages={183--194},
  year={2008},
  publisher={Elsevier},
  doi={https://doi.org/10.1016/j.cie.2007.11.017}
}

@InProceedings{van2018measure,
  author =	{van den Broek, Roel and Hoogeveen, Han and van den Akker, Marjan},
  title =	{{How to Measure the Robustness of Shunting Plans}},
  booktitle =	{18th Workshop on Algorithmic Approaches for Transportation Modelling, Optimization, and Systems (ATMOS 2018)},
  pages =	{3:1--3:13},
  series =	{Open Access Series in Informatics (OASIcs)},
  ISBN =	{978-3-95977-096-5},
  ISSN =	{2190-6807},
  year =	{2018},
  volume =	{65},
  editor =	{Bornd\"{o}rfer, Ralf and Storandt, Sabine},
  publisher =	{Schloss Dagstuhl -- Leibniz-Zentrum f{\"u}r Informatik},
  address =	{Dagstuhl, Germany},
  URL =		{https://drops.dagstuhl.de/entities/document/10.4230/OASIcs.ATMOS.2018.3},
  URN =		{urn:nbn:de:0030-drops-97081},
  doi =		{10.4230/OASIcs.ATMOS.2018.3}
}

@article{jorge1994robustness,
  title={Robustness measures and robust scheduling for job shops},
  author={Jorge Leon, V and David Wu, S and Storer, Robert H},
  journal={IIE transactions},
  volume={26},
  number={5},
  pages={32--43},
  year={1994},
  publisher={Taylor \& Francis},
  doi={https://doi.org/10.1080/07408179408966626}
}

@article{al2005bi,
  title={A bi-objective model for robust resource-constrained project scheduling},
  author={Al-Fawzan, Mohammad A and Haouari, Mohamed},
  journal={International Journal of production economics},
  volume={96},
  number={2},
  pages={175--187},
  year={2005},
  publisher={Elsevier},
  doi={https://doi.org/10.1016/j.ijpe.2004.04.002}
}

@article{kobylanski2007note,
  title={A note on the paper by MA Al-Fawzan and M. Haouari about a bi-objective problem for robust resource-constrained project scheduling},
  author={Kobyla{\'n}ski, Przemys{\l}aw and Kuchta, Dorota},
  journal={International Journal of Production Economics},
  volume={107},
  number={2},
  pages={496--501},
  year={2007},
  publisher={Elsevier},
  doi={https://doi.org/10.1016/j.ijpe.2006.07.012}
}

@article{khemakhem2013efficient,
  title={Efficient robustness measures for the resource-constrained project scheduling problem},
  author={Khemakhem, Mohamed Ali and Chtourou, H{\'e}di},
  journal={International Journal of Industrial and Systems Engineering},
  volume={14},
  number={2},
  pages={245--267},
  year={2013},
  publisher={Inderscience Publishers Ltd},
  doi={https://doi.org/10.1504/IJISE.2013.053738}
}

@article{wilson2014flexibility,
  title={Flexibility and decoupling in simple temporal networks},
  author={Wilson, Michel and Klos, Tomas and Witteveen, Cees and Huisman, Bob},
  journal={Artificial Intelligence},
  volume={214},
  pages={26--44},
  year={2014},
  publisher={Elsevier},
  doi={https://doi.org/10.1016/j.artint.2014.05.003}
}

@article{van2006trade,
  title={The trade-off between stability and makespan in resource-constrained project scheduling},
  author={Van de Vonder, Stijn and Demeulemeester, Erik and Herroelen, Willy and Leus, Roel},
  journal={International Journal of Production Research},
  volume={44},
  number={2},
  pages={215--236},
  year={2006},
  publisher={Taylor \& Francis},
  doi={https://doi.org/10.1080/00207540500140914}
}

@inproceedings{van2013finding,
  title={Finding robust solutions for the stochastic job shop scheduling problem by including simulation in local search},
  author={van den Akker, Marjan and van Blokland, Kevin and Hoogeveen, Han},
  booktitle={Experimental Algorithms: 12th International Symposium, SEA 2013, Rome, Italy, June 5-7, 2013. Proceedings 12},
  pages={402--413},
  year={2013},
  organization={Springer},
  doi={https://doi.org/10.1007/978-3-642-38527-8_35}
}

@article{lambrechts2008tabu,
  title={A tabu search procedure for developing robust predictive project schedules},
  author={Lambrechts, Olivier and Demeulemeester, Erik and Herroelen, Willy},
  journal={International Journal of Production Economics},
  volume={111},
  number={2},
  pages={493--508},
  year={2008},
  publisher={Elsevier},
  doi={https://doi.org/10.1016/j.ijpe.2007.02.003}
}

@article{hazir2010robust,
  title={Robust scheduling and robustness measures for the discrete time/cost trade-off problem},
  author={Haz{\i}r, {\"O}nc{\"u} and Haouari, Mohamed and Erel, Erdal},
  journal={European Journal of Operational Research},
  volume={207},
  number={2},
  pages={633--643},
  year={2010},
  publisher={Elsevier},
  doi={https://doi.org/10.1016/j.ejor.2010.05.046}
}

@article{canon2009evaluation,
  title={Evaluation and optimization of the robustness of dag schedules in heterogeneous environments},
  author={Canon, Louis-Claude and Jeannot, Emmanuel},
  journal={IEEE Transactions on Parallel and Distributed Systems},
  volume={21},
  number={4},
  pages={532--546},
  year={2009},
  publisher={IEEE},
  doi={https://doi.org/10.1109/TPDS.2009.84}
}

@article{boloni2002robust,
  title={Robust scheduling of metaprograms},
  author={B{\"o}l{\"o}ni, Ladislau and Marinescu, Dan C},
  journal={Journal of Scheduling},
  volume={5},
  number={5},
  pages={395--412},
  year={2002},
  publisher={Wiley Online Library},
  doi={https://doi.org/10.1002/jos.115}
}

@inproceedings{shestak2006stochastic,
  title={A stochastic approach to measuring the robustness of resource allocations in distributed systems},
  author={Shestak, Vladimir and Smith, Jay and Siegel, Howard Jay and Maciejewski, Anthony A},
  booktitle={2006 International Conference on Parallel Processing (ICPP'06)},
  pages={459--470},
  year={2006},
  organization={IEEE},
  doi={https://doi.org/10.1109/ICPP.2006.14}
}

@inproceedings{shi2006robust,
  title={Robust task scheduling in non-deterministic heterogeneous computing systems},
  author={Shi, Zhiao and Jeannot, Emmanuel and Dongarra, Jack J},
  booktitle={2006 IEEE international conference on cluster computing},
  pages={1--10},
  year={2006},
  organization={IEEE},
  doi={https://doi.org/10.1109/CLUSTR.2006.311868}
}

@article{nadarajah2008exact,
  title={Exact distribution of the max/min of two Gaussian random variables},
  author={Nadarajah, Saralees and Kotz, Samuel},
  journal={IEEE Transactions on very large scale integration (VLSI) systems},
  volume={16},
  number={2},
  pages={210--212},
  year={2008},
  publisher={IEEE},
  doi={https://doi.org/10.1109/TVLSI.2007.912191}
}

@article{van2005use,
  title={The use of buffers in project management: The trade-off between stability and makespan},
  author={Van de Vonder, Stijn and Demeulemeester, Erik and Herroelen, Willy and Leus, Roel},
  journal={International Journal of production economics},
  volume={97},
  number={2},
  pages={227--240},
  year={2005},
  publisher={Elsevier},
  doi={https://doi.org/10.1016/j.ijpe.2004.08.004}
}

@article{ke2015uncertain,
  title={An uncertain model for RCPSP with solution robustness focusing on logistics project schedule},
  author={Ke, Hua and Wang, Lei and Huang, Hu},
  journal={International Journal of e-Navigation and Maritime Economy},
  volume={3},
  pages={71--83},
  year={2015},
  publisher={Elsevier},
  doi={https://doi.org/10.1016/j.enavi.2015.12.007}
}

@article{fu2015robust,
  title={A robust optimization solution to bottleneck generalized assignment problem under uncertainty},
  author={Fu, Yelin and Sun, Jianshan and Lai, KK and Leung, John WK},
  journal={Annals of Operations Research},
  volume={233},
  pages={123--133},
  year={2015},
  publisher={Springer},
  doi={https://doi.org/10.1007/s10479-014-1631-5}
}

@book{voss2012meta,
  title={Meta-heuristics: Advances and trends in local search paradigms for optimization},
  author={Vo{\ss}, Stefan and Martello, Silvano and Osman, Ibrahim H and Roucairol, Catherine},
  year={2012},
  publisher={Springer Science \& Business Media}
}

@article{xu2013robust,
  title={Robust makespan minimisation in identical parallel machine scheduling problem with interval data},
  author={Xu, Xiaoqing and Cui, Wentian and Lin, Jun and Qian, Yanjun},
  journal={International Journal of Production Research},
  volume={51},
  number={12},
  pages={3532--3548},
  year={2013},
  publisher={Taylor \& Francis},
  doi={https://doi.org/10.1080/00207543.2012.751510}
}

@article{allahverdi2010heuristics,
  title={Heuristics for the two-machine flowshop scheduling problem to minimise makespan with bounded processing times},
  author={Allahverdi, Ali and Aydilek, Harun},
  journal={International Journal of Production Research},
  volume={48},
  number={21},
  pages={6367--6385},
  year={2010},
  publisher={Taylor \& Francis},
  doi={https://doi.org/10.1080/00207540903321657}
}

@article{passage2025new,
  title={A new, efficient approach to speed up local search by estimating the solution quality: an application to stochastic, parallel machine scheduling: G. Passage et al.},
  author={Passage, Guido and van den Akker, Marjan and Hoogeveen, Han},
  journal={Journal of Heuristics},
  volume={31},
  number={3},
  pages={26},
  year={2025},
  publisher={Springer},
  doi={https://doi.org/10.1007/s10732-025-09562-5}
}

@article{liang2020integration,
  title={The integration of resource allocation and time buffering for bi-objective robust project scheduling},
  author={Liang, Yangyang and Cui, Nanfang and Hu, Xuejun and Demeulemeester, Erik},
  journal={International Journal of Production Research},
  volume={58},
  number={13},
  pages={3839--3854},
  year={2020},
  publisher={Taylor \& Francis},
  doi={https://doi.org/10.1080/00207543.2019.1636319}
}

@article{ma2022genetic,
  title={A genetic algorithm for proactive project scheduling with resource transfer times},
  author={Ma, Zhiqiang and Zheng, Weibo and He, Zhengwen and Wang, Nengmin and Hu, Xuejun},
  journal={Computers \& Industrial Engineering},
  volume={174},
  pages={108754},
  year={2022},
  publisher={Elsevier},
  doi={https://doi.org/10.1016/j.cie.2022.108754}
}

@article{lambrechts2011time,
  title={Time slack-based techniques for robust project scheduling subject to resource uncertainty},
  author={Lambrechts, Olivier and Demeulemeester, Erik and Herroelen, Willy},
  journal={Annals of operations Research},
  volume={186},
  number={1},
  pages={443--464},
  year={2011},
  publisher={Springer},
  doi={doi.org/10.1007/s10479-010-0777-z}
}

@article{lambrechts2008proactive,
  title={Proactive and reactive strategies for resource-constrained project scheduling with uncertain resource availabilities},
  author={Lambrechts, Olivier and Demeulemeester, Erik and Herroelen, Willy},
  journal={Journal of scheduling},
  volume={11},
  number={2},
  pages={121--136},
  year={2008},
  publisher={Springer},
  doi={https://doi.org/10.1007/s10951-007-0021-0}
}

@article{ma2019computational,
  title={A computational experiment to explore better robustness measures for project scheduling under two types of uncertain environments},
  author={Ma, Zhiqiang and Demeulemeester, Erik and He, Zhengwen and Wang, Nengmin},
  journal={Computers \& Industrial Engineering},
  volume={131},
  pages={382--390},
  year={2019},
  publisher={Elsevier},
  doi={https://doi.org/10.1016/j.cie.2019.04.014}
}

@article{xiao2026branch,
  title={A branch-and-price method for parallel machine scheduling with stochastic processing times},
  author={Xiao, Fan and Gao, Yangyang and Liang, Zhe},
  journal={Optimization Letters},
  pages={1--39},
  year={2026},
  publisher={Springer},
  doi={https://doi.org/10.1007/s11590-026-02301-0}
}

@article{bruni2011heuristic,
  title={A heuristic approach for resource constrained project scheduling with uncertain activity durations},
  author={Bruni, Maria Elena and Beraldi, Patrizia and Guerriero, Francesca and Pinto, Erika},
  journal={Computers \& Operations Research},
  volume={38},
  number={9},
  pages={1305--1318},
  year={2011},
  publisher={Elsevier},
  doi={https://doi.org/10.1016/j.cor.2010.12.004}
}

@article{liu2021parallel,
  title={Parallel machine scheduling with stochastic release times and processing times},
  author={Liu, Xin and Chu, Feng and Zheng, Feifeng and Chu, Chengbin and Liu, Ming},
  journal={International Journal of Production Research},
  volume={59},
  number={20},
  pages={6327--6346},
  year={2021},
  publisher={Taylor \& Francis},
  doi={https://doi.org/10.1080/00207543.2020.1812752}
}

@article{skutella2005stochastic,
  title={Stochastic machine scheduling with precedence constraints},
  author={Skutella, Martin and Uetz, Marc},
  journal={SIAM Journal on Computing},
  volume={34},
  number={4},
  pages={788--802},
  year={2005},
  publisher={SIAM},
  doi={https://doi.org/10.1137/S0097539702415007}
}

\end{document}